\documentclass[reqno,english,11pt]{amsart}
\usepackage[]{amsmath,amssymb,amsfonts,latexsym,amsthm,enumerate}
\usepackage{amssymb,bbm,mathrsfs,tikz,stmaryrd}
\usetikzlibrary{calc,decorations.pathmorphing}
\usetikzlibrary{patterns}
\usepackage[normalem]{ulem}
\usepackage{enumerate,graphicx,paralist}
\usepackage{upgreek}
\usepackage[font=small]{caption}
\usepackage[noadjust]{cite}
\usepackage[perpage]{footmisc}

\numberwithin{equation}{section}
\usepackage[psdextra=true,colorlinks=true, pdfstartview=FitV, linkcolor=blue,
  citecolor=blue, urlcolor=blue,pagebackref=false]{hyperref}
\usepackage[shortlabels]{enumitem}  
\usepackage[final]{showkeys}

\usepackage{xspace}
\usepackage[letterpaper]{geometry}
\usepackage[colorinlistoftodos]{todonotes}
\presetkeys{todonotes}{inline, color=green}{}
\usepackage{comment}
\usepackage[capitalise]{cleveref}
\newtheorem{maintheorem}{Theorem}  
\newtheorem{maincoro}[maintheorem]{Corollary}

\newtheorem{theorem}{Theorem}[section]
\newtheorem*{theorem*}{Theorem}
\newtheorem{lemma}[theorem]{Lemma}
\newtheorem{claim}[theorem]{Claim}
\newtheorem{proposition}[theorem]{Proposition}
\newtheorem{observation}[theorem]{Observation}
\newtheorem{fact}[theorem]{Fact}

\newtheorem{remark}[theorem]{Remark}
\newtheorem*{remark*}{Remark}

\theoremstyle{definition}{
\newtheorem{example}[theorem]{Example}
\newtheorem{definition}[theorem]{Definition}
\newtheorem*{definition*}{Definition}

\newtheorem*{question*}{Question}
\newtheorem*{example*}{Example}
\newtheorem*{examples*}{Examples}
}

\newcommand{\E}{\mathbb E}

\newcommand{\R}{\mathbb R}

\newcommand{\Z}{\mathbb Z}

\newcommand{\cD}{\mathcal D}

\newcommand{\cA}{{\mathcal A}}

\newcommand{\cC}{{\mathcal C}}
\newcommand{\cE}{{\mathcal E}}
\newcommand{\cF}{{\mathcal F}}

\newcommand{\cH}{{\mathcal H}}
\newcommand{\cI}{{\mathcal I}}

\newcommand{\cP}{{\mathcal P}}

\newcommand{\cS}{{\mathcal S}}

\newcommand{\cX}{{\mathcal X}}

\newcommand{\fb}{\mathfrak b}

\newcommand{\fg}{\mathfrak g}
\newcommand{\fm}{\mathfrak m}

\newcommand{\fC}{\mathfrak C}

\newcommand{\fP}{\mathfrak P}

\newcommand{\fX}{\mathfrak X}

\newcommand{\sw}{\mathsf w}

\newcommand{\sC}{\mathsf C}

\newcommand{\scC}{\mathscr C}

\newcommand{\sign}{\operatorname{sgn}}
\newcommand{\Cov}{\operatorname{Cov}}

\DeclareMathOperator{\diam}{diam}
\DeclareMathOperator{\dist}{dist}
\DeclareMathOperator{\supp}{supp}

\newcommand{\one}{\mathbbm{1}}

\renewcommand{\restriction}{\mathord{\upharpoonright}}
\renewcommand{\epsilon}{\varepsilon}

\DeclareMathOperator{\Int}{Int}
\DeclareMathOperator{\Ext}{Ext}

\newcommand{\starsim}{\overset{*}{\sim}}
\newcommand{\agreeoutside}[1]{\overset{\;#1}{\equiv}}
\newcommand{\partialvtxin}{\partial^{\mathrm{in}}_{\mathtt v}}
\newcommand{\partialvtxout}{\partial^{\mathrm{ex}}_{\mathtt v}}
\newcommand{\partialedge}{\partial_{\mathtt e}}
\newcommand{\partialhat}{\widehat\partial_{\mathtt v}}
\newcommand{\fp}{\mathfrak p}
\newcommand{\scS}{\mathscr S}
\newcommand{\Cst}{C_\star}
\newcommand{\shadow}{\Delta}
\newcommand{\scV}{\mathscr V}

\newcommand{\given}{\;\big|\;}

\makeatletter
\crefname{maintheorem}{Theorem}{Theorems}
\crefname{maincoro}{Corollary}{Corollaries}
\AddToHook{env/maincoro/begin}{\crefalias{maintheorem}{maincoro}}
\crefname{step}{Step}{Steps}
\crefname{part}{Part}{Parts}
\crefname{case}{Case}{Cases}
\crefname{claim}{Claim}{Claims}
\AddToHook{env/lemma/begin}{\crefalias{theorem}{lemma}}
\AddToHook{env/claim/begin}{\crefalias{theorem}{claim}}
\AddToHook{env/proposition/begin}{\crefalias{theorem}{proposition}}
\AddToHook{env/corollary/begin}{\crefalias{theorem}{corollary}}
\AddToHook{env/definition/begin}{\crefalias{theorem}{definition}}
\AddToHook{env/example/begin}{\crefalias{theorem}{example}}
\AddToHook{env/question/begin}{\crefalias{theorem}{question}}
\AddToHook{env/remark/begin}{\crefalias{theorem}{remark}}
\AddToHook{env/observation/begin}{\crefalias{theorem}{observation}}
\AddToHook{env/fact/begin}{\crefalias{theorem}{fact}}
\AddToHook{cmd/appendix/before}{%
  \crefalias{section}{appendix}%
  \crefalias{subsection}{appendix}%
}
\makeatother

\usepackage{crossreftools}
\title[Cluster expansion for long-range tilts]{Radius decay and boundary polymers: \\ cluster expansion for long-range many-body tilts}

\author{Beno\^{\i}t Laslier}
\address{B.\ Laslier\hfill\break
Universit\'e de Rouen Normandie\\
LMRS\\
76801 Saint-\'Etienne-du-Rouvray\\ France.}
\email{benoit.laslier@univ-rouen.fr}
\author{Eyal Lubetzky}
\address{E.\ Lubetzky\hfill\break
Courant Institute\\ New York University\\
251 Mercer Street\\ New York, NY 10012, USA.}
\email{eyal@courant.nyu.edu}

\begin{document}

\begin{abstract}
We consider a nearest-neighbor Gibbs measure on a region in $\Z^2$ (e.g., the Ising model) tilted by $\sum_{x,r} g_{x,r}$, where $g_{x,r}$ is an arbitrary function of the spins in the ball of radius $r$ about $x$, with $\|g_{x,r}\|_\infty\leq \exp(-C r)$: an exponential decay in the radius of the ball rather than in its volume. Standard cluster expansions fail for such many-body tilts, as a ball of radius $r$ carries $\exp(c r^2)$ configurations. We show that the partition function can be resummed into a polymer model whose polymers are outer boundaries of connected unions of balls, decorated by the spins along them, with weights decaying exponentially in the boundary length.
The proof uses only the domain Markov property of the base measure; it develops a version of this resummation argument that appeared in a recent work of the authors in the context of RG analysis of Grassmann variables.
As~applications, we show convergent cluster expansions for many-body tilts of Ising models at high and low temperatures; existence of infinite-volume limits; and exponential decay of correlations.
\end{abstract}

\maketitle
 \vspace{-0.6cm}

\section{Introduction}\label{sec:intro}

Consider a nearest-neighbor Gibbs measure $\mu_\Lambda$ on a finite region $\Lambda\subset\Z^2$, assigning a spin $\sigma_x\in \scS$ to each site $x\in \Lambda$, whose Hamiltonian $\cH_\Lambda:\scS^\Lambda\to\R$ features generic (potentially inhomogeneuous) potential on the edges, external fields and boundary conditions (formally given in \S\ref{subsec:setting}, \cref{eq:gibbs-def}):
\begin{equation}\label{eq:base-def}
  \mu_\Lambda(\sigma) \propto \exp\left(-\cH_\Lambda(\sigma)\right)\qquad\mbox{for }\sigma\in\scS^\Lambda\,.
\end{equation} 
The goal is to study a perturbation of $\mu_\Lambda$ via long-range interaction functions $g_{x,r}:\scS^{B(x,r)}\to\R$ 
for all $x\in\Lambda$ and $r\geq 0$, where $B(x,r)=\{y\in \Z^2:\|y-x\|_\infty\leq r\}$ is the $\ell_\infty$-ball of radius $r$ about $x$. In particular, we aim to establish a convergent cluster expansion formula for the partition function of the tilted model when $ | g_{x, r} | \le e^{-\Cst (r+1)}$. Write $\fg=(g_{x,r})_{x\in\Lambda,\, r\geq 0}$  and set
\begin{equation}\label{eq:HD-def}
\cH_\Lambda^\fg(\sigma) := \sum_{(x,r)\,:\;B(x,r)\subset \Lambda} g_{x,r}\big(\sigma\restriction_{B(x,r)}\big)\,,\qquad
Z_\Lambda^{\fg} := \E\big[e^{\cH_\Lambda^\fg(\sigma)}\big]\,.
\end{equation} The object of interest is the tilted measure
\begin{equation}\label{eq:mu-def}
\mu_\Lambda^{\fg}(\sigma) := \frac{\exp(\cH^\fg_\Lambda(\sigma))}{Z_\Lambda^{\fg}}\,\mu_\Lambda(\sigma)\qquad\mbox{for $\sigma\in \scS^\Lambda$}\,.
\end{equation}
Our standing assumption on the interaction term $\fg$ is that, for some $\Cst>0$,
\begin{equation}\label{eq:g-decay}
\|g_{x,r}\|_\infty\leq e^{-\Cst (r+1)}\qquad\mbox{for all $x\in\Lambda$ and $r\geq 0$}\,.
\end{equation}
Cluster expansion refers to the following formal equality for the partition function of a model of abstract polymers $\gamma \in \cP$ with weights $\sw(\gamma)$ and interactions $\zeta(\gamma,\gamma')$:
\begin{equation}\label{eq:cluster_expansion}
 \sum_{ \{ \gamma_i \} \subset \Gamma} \prod_i \sw(\gamma_i) \prod_{i, j} \zeta(\gamma_i, \gamma_j) = \exp\bigg( \sum_{n\geq1}\frac1{n!} \sum_{(\gamma_1, \ldots, \gamma_n)} \phi_\zeta( \gamma_1, \ldots, \gamma_n)\prod_{i} \sw(\gamma_i)\bigg)\,,
\end{equation}
where the $\phi_\zeta$ are Ursell functions: explicit functions that account for the nontrivial combinatorial terms arising from the pairwise interactions (see \cref{eq:ursell-def}); a version of this goes back to Ursell~\cite{Ursell27} and the work of Mayer~\cite{Mayer37} on the virial expansion; Gruber and Kunz~\cite{GruberKunz71} introduced polymers; the abstract form used here is due to Kotecký and Preiss \cite{KoteckyPreiss86}. See the comprehensive account~\cite[\S5]{FriedliVelenik18}. 

The equality is only formal in some settings because the series on the right-hand side may be divergent even when the set of polymers $\Gamma$ is finite, but when it holds---typically at either low or high temperature---it is quite useful (see, e.g., \cite{MinlosSinai67,Dobrushin72,PirogovSinai75,PirogovSinai76} for early applications of this technique). Indeed, the right-hand side gives an explicit expression of the free energy from which it is easy to derive important properties of the associated measure, such as analyticity of the free energy or local statistics. Over the long history of study of this topic, many criteria were introduced for the convergence of the expansion (see \cite{FernandezProcacci07} for a comparison of multiple such criteria, including \cite{Dobrushin96,KoteckyPreiss86}); the vanilla criterion (see \cite[\S5.4]{FriedliVelenik18}) is to find a function $a(\cdot)$ on polymers such that, for every $\gamma_0\in\Gamma$,
\begin{equation}\label{eq:textbook_bound}
\sum_{\gamma\in\Gamma} |\sw(\gamma)|\, e^{a(\gamma)}\, |\zeta( \gamma, \gamma_0) -1|\leq a( \gamma_0)\,.
\end{equation}

The simplest application of \cref{eq:cluster_expansion} to our setting is when the base measure $\mu_\Lambda$ is a product measure: 
Writing $e^{g_{x, r}} =1 + (e^{g_{x, r}} -1)$ and developing the sum over $\{(x,r): B(x,r)\subset\Lambda\} $ leads to 
\begin{equation}\label{eq:classical-expansion}
Z_\Lambda^{\fg} = \sum_{ \{ (x_i, r_i) \}} \E \bigg[\prod_i \left(e^{ g_{x_i, r_i} ( \sigma \restriction_{B(x_i, r_i)})} -1\right)\bigg]\,.
\end{equation}
Using the independence in the base measure $\mu_\Lambda$, this can easily be rewritten as in \cref{eq:cluster_expansion} with polymers $\gamma$ given by connected unions of balls and the pairwise interaction $\zeta( \gamma_1, \gamma_2) = \one_{\{\gamma_1 \cap \gamma_2 = \emptyset\}}$. The convergence criterion~\labelcref{eq:textbook_bound} is not hard to check if $|g_{x, r}| \leq e^{-C r^2}$ for a large enough $C>0$ but can fail with the weaker $e^{-\Cst r}$ (radius as opposed to area decay); see \cref{rem:naive} for more on this. 

An example to keep in mind is the following: consider a high temperature Ising model, tilted by $\cH_\Lambda^\fg(\sigma)= \sum_x \exp[-C \diam(\cC_x(\sigma))]$, where $\diam(\cC_x)$ is the diameter of the connected cluster of same-spin sites that contains $x$. The decay rate of $\exp(-\diam(\cC_x))$ is insufficient to combat the $\exp( c |\cC_x|)$ number of possible connected component in the classical expansion mentioned above.

Our main theorem provides an alternative way to rewrite $Z_\Lambda^\fg$  as a sum of polymers which does satisfy~\labelcref{eq:textbook_bound} notwithstanding the many-body interaction~$\fg$ that decays exponentially in the radius rather than area. Throughout, two sites $x,y$ are adjacent, denoted $x\sim y$, if $\|x-y\|_1=1$, and $*$-adjacent, denoted $x\starsim y$, if $\|x-y\|_\infty=1$, extending to connectivity and $*$-connectivity of sets and the diameters $\diam(\cdot)$ and $\diam_*(\cdot)$; a finite $\Lambda\subset \Z^2$ is simply connected if $\Z^2\setminus\Lambda$ is $*$-connected. 

\begin{maintheorem}\label{thm:main}
There exists an absolute constant $C_0>0$ such that the following holds. Let $\mu_\Lambda$ be a generic nearest-neighbor Gibbs measure (as per \cref{eq:gibbs-def}) on a finite simply connected $\Lambda\subset \Z^2$. Let $\fg$ be a family of interaction functions satisfying \cref{eq:g-decay} for a fixed $\Cst\geq C_0$.  Then there exist a family $\Gamma$ of connected subsets of $\Lambda$, a symmetric \emph{compatibility} relation $\perp$ on $\Gamma$ such that $\gamma_1\perp\gamma_2$ implies $\gamma_1\cap\gamma_2=\emptyset$ and $\gamma_1\not\perp\gamma_2$ implies the graph distance in $\Z^2$ ($\ell_1$ distance) of $\gamma_1,\gamma_2$ is at most~$1$, and a weight function $\sw:\cP\to\R$ on the polymers $\cP=\{(\gamma,\eta):\, \gamma\in\Gamma,\, \eta\in \scS^\gamma\}$, such that 
\begin{equation}\label{eq:polymer-rep}
Z_\Lambda^{\fg} = \sum_{\cF}\E\bigg[\prod_{(\gamma,\eta)\in\cF} \sw(\gamma,\eta)\,\one_{\{\sigma\restriction_\gamma = \eta\}}\bigg]\,,
\end{equation}
where the sum runs over all families $\cF=\{(\gamma_i, \eta_i)\}\subset\cP$ of polymers with $\gamma_i \perp \gamma_j $ for all $i\neq j$ (including the empty family), and
\begin{equation}\label{eq:w-bound}
|\sw(\gamma,\eta)|\leq e^{-\Cst|\gamma|/10}\qquad\mbox{for every polymer $(\gamma,\eta)$}\,.
\end{equation}
\end{maintheorem}

The definitions of the polymers $\cP$ and weights $\sw$ are completely constructive, but do require special care, e.g., in order to maintain the domain Markov property that our recursive expansion requires, as well as to allow nonnegative weights $\sw(\cdot)$.
We thus postpone these definitions to the main proof, and at this point just note that $\Gamma$, the set of \emph{polymer supports}, does not depend on~$\mu_\Lambda$. The main idea is that a polymer $(\gamma, \eta)$ will describe the boundary of a connected set of balls instead of the full domain, so that its size will be controlled by the sum of the radii of the involved balls. The weight of a polymer is then essentially the contribution of all the boundary balls, formally expressed as a ratio of two partition functions. The idea of using polymer boundaries with weights given by ratios of interior partition functions goes back to Pirogov--Sinai theory in the form given to it by Zahradn\'{\i}k~\cite{Zahradnik84} in 1984. While we use the same overarching recipe here, our procedure captures different applications (which neither contain nor are contained in those of Pirogov--Sinai).

Very little effort was made to optimize the constant $C_0$; the current proof is valid for $C_0=35$.
The simplest application of the theorem is for the case where the reference measure $\mu_\Lambda$ is i.i.d., whence the expectation in \cref{eq:polymer-rep} factorizes and gives a convergent cluster expansion formula for~$Z_\Lambda^\fg$. Note that this case covers arbitrary finite-range spin system at high enough temperatures, as their interaction terms can be embedded into the functions $g_{x,r}$ in $\fg$. However, when the base measure is not i.i.d., one can still use the cluster expansion formula for all given spin configurations to get a rewriting of the partition function with a ``nicer'' effective potential.
\begin{maincoro}\label{cor:potential}
In the setting of \cref{thm:main}, let $\delta:=\Cst/20-2$. There exist functions $\Psi_A:\scS^A\to\R$, indexed by connected sets $A\subset\Lambda$, such that
\begin{equation}\label{eq:effective-potential}
Z_\Lambda^{\fg} = \E\bigg[\exp\Big(\sum_{A}\Psi_A(\sigma\restriction_A)\Big)\bigg]
\quad\mbox{and}\quad
\sum_{A\ni y}\|\Psi_A\|_\infty\, e^{\delta|A|}\leq e^{-\delta}\ \ \mbox{for every $y\in\Lambda$}\,.
\end{equation}
\end{maincoro}
A version of \cref{thm:main} specialized to dimer configurations, as well as a version of \cref{cor:potential} specialized to its associated system of Grassmann variables, appeared in a recent work of the authors~\cite{LaslierLubetzky24} (see Prop.~8.11 there). Here we prove a more generalized, but also slightly improved, variant of the argument, dispensing with the restriction that the functions $g_{x, r}$ must be nonnegative. Note that in that context, the base measure (uniform dimers) is critical, so the fact that there are no stability or correlation assumptions on $\mu_\Lambda$ is crucial. 

We stress that \cref{eq:effective-potential} is foremost an identity for the partition function $Z_\Lambda^\fg$: while it implies that  $\mu_\Lambda^\fg\propto e^{\cH^\fg_\Lambda}\mu_\Lambda$ and $\tilde\mu_\Lambda^\fg \propto e^{\sum_A\Psi_A}\mu_\Lambda$ have the same normalizer $Z_\Lambda^\fg$, these are two different measures.
Nevertheless, one can recover local information on the base measure by using the standard recipe: for an event $E$ measurable with respect to $ B(x_0, r_0)$ for some $x_0,r_0$, modifying $\tilde g_{x_0,r_0} := g_{x_0,r_0} + \lambda \one_E$ and letting $Z^{\tilde\fg}_\Lambda$ be the associated partition function (all other $g_{x,r}$ unchanged), one has
\[
\mu_\Lambda^{\fg}(E) = \frac{\partial}{\partial \lambda} \log Z^{\tilde\fg}_\Lambda\big|_{\lambda=0}\,.
\]
\begin{remark}[Analyticity or lack thereof]\label{rem:analyticity}
 While classical cluster expansion readily implies the analyticity of the free energy with respect to a global prefactor for the interactions (temperature), surprisingly, we do not expect such analyticity to hold in general under \cref{eq:g-decay}...
 This is because the simple scaling $\fg \to \beta \fg$ of the interactions has a non-trivial effect on the resulting weights $\sw$ and $\partial_\beta\sw(\gamma)$ may diverge with $|\gamma|$. We will see in \cref{rem:analyticity-obstruction} a detailed heuristic based on the cumulants showing that the free energy should not be analytic in general. 

 However, we will still prove that one can recover local events with the technique given above since a local perturbation of the weights such as $\tilde g_{x_0,r_0} := g_{x_0,r_0} + \lambda \one_E$ creates an analytic perturbation of the free energy: $\log Z^{\tilde\fg}_\Lambda  - \log Z^{\fg}_\Lambda$ will converge as $\Lambda \to \infty$ to an analytic function of~$\lambda$.
\end{remark}

A recent improvement on cluster expansion for the Ising model with long-range interactions on~$\Z^d$ for arbitrary $d$ was obtained by Affonso et al.\ in the work~\cite{ABEH25}, and its sequel~\cite{ABMRW25}. These works can cover stronger interactions (decaying only as polynomials of the distance), while we are restricted to exponential decay and to $\Z^d$ for $d=2$ only; however, the works~\cite{ABEH25,ABMRW25} can only handle a two-body interaction term, not a general function (arbitrary many-body interaction).

In fact, we argue that in a sense, one should not think of our result as addressing long-range interactions but rather as addressing ``many many''-body interactions,  quite orthogonally to \cite{ABMRW25,ABEH25}. Indeed, consider a model of lattice polymer where every set $S$ is given a weight $\sw(S)$ (which one could think of as coming from a naive expansion of a system with many-body interactions), e.g., \[\sw(S) := \frac{\exp(-C_0 |S|) }{\min_{T} \prod_{(xy) \in T} |x - y|^{4}}\,,\] where the minimum over $T$ runs over all abstract trees with vertices indexed by all points of $S$.
For a large enough $C_0>0$, this choice satisfies the criterion in~\labelcref{eq:textbook_bound}, despite the polynomially decaying interactions. In other words, an interaction decaying as $|x-y|^{-4}$ is not problematic at all as long as the weight decays fast enough, exponentially, with the number of vertices.
Compared to that case, the naive expansion in our model ($\mu_\Lambda^\fg$ from \cref{eq:mu-def}) involves sets $S$ that are connected, but their weights only decay as $\exp(-C_0 \sqrt{| S|})$ which is not enough for the classical criterion~\labelcref{eq:textbook_bound}.

\begin{remark}[The direct expansion]\label{rem:naive}
    The classical criterion \labelcref{eq:textbook_bound} may fail for the natural choice of using connected unions of balls as polymers. However, this still does not imply that the cluster expansion is divergent for these polymers (more involved convergence criteria exist in the literature). See \cref{prop:naive} for a scenario where the actual expansion for these polymers diverges.
\end{remark}

\subsection*{Applications}

The next examples set the domain $\Lambda$ to be an $N_1\times N_2$ box in $\Z^2$. Recall that the Ising model on $\Lambda$ at inverse temperature $\beta\geq0$ is the probability measure on $\{\pm1\}^\Lambda$ given by
\begin{equation}\label{eq:ising-def}
\mu_\Lambda(\sigma)\propto\exp\Big(-\beta\sum_{xy \in E(\Lambda)}\one_{\{\sigma_x\neq\sigma_y\}}\Big)\,,
\end{equation}
where $E(\Lambda)$ is the edge set of $\Lambda$. Given a boundary condition $\tau$, we further sum over edges $xy$ with $x\in\Lambda$ and $y\notin\Lambda$, where $\sigma_y$ is prescribed to be $\tau_y$. (See also \cref{eq:gibbs-def} with $\scV_{x,y}(s_1,s_2)=\beta\one_{\{s_1\neq s_2\}}$.)

\begin{example}[High temperature Ising]\label{ex:ising-high}
Let $\mu_\Lambda$ be the Ising model of \cref{eq:ising-def} with an arbitrary boundary condition. Let $\scC_x=\scC_x(\sigma)$ be the connected component of $x\in\Lambda$ in $\{y\in\Lambda:\sigma_y=\sigma_x\}$, and, for a large enough constant $C>0$, consider
\begin{equation}\label{eq:ising-tilt-cluster-diam}
\cH_\Lambda^\fg(\sigma) := \sum_{x\in\Lambda} e^{-C(\diam(\scC_x)+1)} 
\end{equation}
(equivalently, $\cH_\Lambda^\fg(\sigma)=\sum_{\scC} |\scC| e^{-C(\diam(\scC)+1)}$, the sum running over the connected components $\scC$ of the sets of plus spins and of minus spins), as per \cref{eq:HD-def,eq:g-decay}. If $\beta$ is small enough, then $\log Z_\Lambda^\fg$ is given by absolutely convergent cluster expansions, uniformly in $\Lambda$. Moreover, the tilted measures $\mu_\Lambda^\fg$ converge weakly, as $\Lambda\uparrow\Z^2$, to a translation-invariant probability measure~$\mu_\infty^\fg$ on $\{\pm1\}^{\Z^2}$ (independent of the boundary conditions), which has exponential decay of correlations: there are absolute constants $C',c>0$ so that, for all finite disjoint subsets $X,Y\subset\Z^2$,
\begin{equation}\label{eq:ex-decay}
\sup_{\xi,\xi'}\Big|\mu_\infty^\fg(\sigma\restriction_{X}=\xi,\ \sigma\restriction_{Y}=\xi')-\mu_\infty^\fg(\sigma\restriction_{X}=\xi)\,\mu_\infty^\fg(\sigma\restriction_{Y}=\xi')\Big|\leq e^{C'(|X|+|Y|)}\,e^{-c\,\dist(X,Y)}\,.
\end{equation}
\end{example}

\begin{example}[Low temperature Ising]\label{ex:ising-low}
In the setting of \cref{ex:ising-high} (Ising tilted by \labelcref{eq:ising-tilt-cluster-diam}), consider instead all-plus boundary conditions. If $\beta$ is large enough, then the same conclusions hold: $\log Z_\Lambda^\fg$ is given, up to an explicit sum of local terms, by absolutely convergent cluster expansions, uniformly in $\Lambda$, and the tilted measures $\mu_\Lambda^\fg$ converge weakly as $\Lambda\uparrow\Z^2$ to a translation-invariant probability measure satisfying \cref{eq:ex-decay}.
\end{example}

\begin{example}[Uniform measure]\label{ex:product-abs-mag}
Let $\mu_\Lambda$ be the uniform measure on $\{\pm1\}^\Lambda$ (equivalent to the Ising model at $\beta=0$), and, for a constant $C>0$, define $\mu_\Lambda^\fg$ via the tilt $\cH_\Lambda^\fg$ given by 
\[
g_{x,r}(\sigma):=e^{-C(r+1)}\,\Big|\sum_{v\in B(x,r)}\sigma_v\Big|\qquad\mbox{for all $x$ and $r\geq 0$}\,.
\]
Then, if $C$ is large enough, the conclusions of \cref{ex:ising-high} hold true: $\log Z_\Lambda^\fg$ has an absolutely convergent cluster expansion, the measures $\mu_\Lambda^\fg$ converge to an infinite-volume limit, which has exponential decay of correlations. However, the classical cluster expansion for this model---polymers are connected unions of balls---\emph{fails} to converge: the sum over $n\geq 1$ in \cref{eq:cluster_expansion} diverges already in a finite box (driven by a polymer of weight $>1$ repeated in the ordered $n$-tuples); see \cref{prop:naive}.
\end{example}

\Cref{ex:ising-high,ex:ising-low} are special cases of \cref{prop:product,prop:lowT} in \cref{sec:corollaries}, respectively. The former treats product base measures, and extends to any nearest-neighbor model at sufficiently high temperature by embedding the nearest-neighbor interaction in the tilt at $r=1$ (see \cref{rem:highT}); in particular, it applies to \cref{ex:ising-high,ex:product-abs-mag}. The latter treats the Ising model at low temperature under all-plus boundary conditions and an arbitrary spin-flip symmetric tilt $\fg$.

\section{Proof of \texorpdfstring{\cref{thm:main}}{Theorem 1} and \texorpdfstring{\cref{cor:potential}}{Corollary 2}}\label{sec:boundary-polymers}

\subsection{Setup}\label{subsec:setting}
Throughout this paper, we let $\dist(\cdot,\cdot)$ be the graph distance in $\Z^2$ ($\ell_1$ distance), and let $\dist_*(\cdot,\cdot)$ be the graph distance with respect to $*$-adjacency ($\ell_\infty$ distance). For a finite $V\subset\Z^2$, let $E(V) := \{xy\,:\; x,y\in V \mbox{ with } x\sim y\}$ denote the edges internal to~$V$, denote its boundary edges by $\partialedge V := \{xy\,:\; x\in V,\, y\in V^c,\, x\sim y\} $, and let $\partialvtxin V := \{x\in V:\; xy\in \partialedge V\mbox{ for some }y\}$ and $\partialvtxout V:=\{y\in V^c:\; xy\in\partialedge V \mbox{ for some }x\}$ denote its internal and external vertex boundaries. Recall that $B(x,r)=\{y\in\Z^2:\dist_*(x,y)\leq r\}$ for $x\in\Z^2$ is a box with $(2r+1)^2$ sites, of which exactly $(2r+1)^2-(2r-1)^2=8r$ sites are $*$-adjacent to some $z\notin B(x,r)$ when $r\geq 1$.

Fix a finite set $\scS$ of spin values and a finite region $V\subset \Z^2$, and consider a family $\scV=(\scV_{x,y},\scV_x)$ of potential functions: symmetric $\scV_{x,y}:\scS\times \scS\to\R\cup\{+\infty\}$ as well as $\scV_x:\scS\to\R\cup\{+\infty\}$ (where the $+\infty$ value corresponds to a hard constraint), indexed by pairs $x,y$ corresponding to edges in $E(V)\cup \partialedge V$ (ones incident to $V$) and by sites $x\in V$, respectively. The associated Gibbs measure on $V$ with boundary condition $\tau\in \scS^{V^c}$ is the probability measure on $\scS^V$ given by
\begin{equation}\label{eq:gibbs-def}
\mu_V^\tau(\sigma) := \frac{e^{-\cH_V^\tau(\sigma)}}{Z_V^\tau} \quad\mbox{for}\quad \cH_V^\tau(\sigma):=\sum_{xy\in E(V) \cup \partialedge V}\scV_{x,y}(\sigma_x,\sigma_y)+\sum_{x\in V}\scV_x(\sigma_x)\,,
\end{equation}
where we only consider $\sigma$ that agree with $\tau$ on $V^c$, and $Z_V^\tau = \sum_{\sigma\in\scS^V} e^{-\cH^\tau_V(\sigma)}$ is the normalizing constant, assumed to be nonzero (there exists a legal $\sigma$ with respect to the hard constraints); we write $\E_V^\tau$ for the corresponding expectation. Note that $\cH_V^\tau$, $Z_V^\tau$ and $\mu_V^\tau$ depend on $\tau$ only through $\tau\restriction_{\partialvtxout V}$. Accordingly, if $\xi\in\scS^W$ for some $W$ with $\partialvtxout V\subset W\subset V^c$, we write $\cH_V^\xi$, $Z_V^\xi$ and $\mu_V^\xi$ for $\cH_V^\tau$, $Z_V^\tau$ and $\mu_V^\tau$ with any $\tau\in\scS^{V^c}$ extending~$\xi$. We call the boundary condition \emph{free} when $\scV_{x,z}\equiv0$ for every edge $xz\in\partialedge V$. Finally, for $V\subset\Lambda$ we define $\cH_V^\fg$ and $Z_V^{\tau,\fg} := \E_V^\tau[e^{\cH_{V}^\fg(\sigma)}]$ as in \cref{eq:HD-def} when replacing $\Lambda$ by $V$, and $\mu_V^{\tau,\fg}$ as the tilt of $\mu_V^\tau$ via $\cH_{V}^\fg$ as in \cref{eq:mu-def}. Once the boundary condition $\tau$ is fixed, we may drop its superscript from the notation for brevity. 

The principal examples are the Ising and Potts measures with free or fixed boundary conditions, for which $\scV_{x,y}(s,s')=-\beta\one_{\{s=s'\}}$, and wherein the $\beta=0$ case is the product measure ($\scV_{x,y}\equiv0$).

Our proof will hinge on the well-known Domain Markov property of Gibbs measures, stated in the following lemma. Its proof is routine and we provide it for completeness in \cref{app:markov}.
\begin{lemma}[Domain Markov property]\label{lem:markov}
Let $V\subset\Z^2$ be finite and let $\tau\in \scS^{V^c}$ satisfy $Z_V^\tau>0$. Let $U\subset V$, let $\xi\in \scS^{U^c}$ satisfy $\xi\restriction_{V^c}=\tau$, and suppose that $\mu_V^\tau(\cE_\xi)>0$, where $\cE_\xi$ is the event that $\sigma\restriction_{V\setminus U} = \xi\restriction_{V\setminus U}$. Then $Z_U^{\xi}>0$ and, for every $F:\scS^U\to\R$,
\[
\E_V^\tau\big[F(\sigma\restriction_U)\,\one_{\cE_\xi}\big] = \mu_V^\tau(\cE_\xi)\,\E_U^{\xi}[F]\,.
\]
Moreover, if $U$ is partitioned into sets $U_1,\ldots,U_m$ with $\dist(U_i,U_j)\geq 2$ for all $i\neq j$, then $Z_{U_k}^{\xi}>0$ for every $k$ and $\mu_U^{\xi}=\bigotimes_{k=1}^m\mu_{U_k}^{\xi}$.
\end{lemma}

\subsection{Geometry of boundaries}\label{subsec:geometry}

For a finite set $K\subset\Z^2$, the set $\Z^2\setminus K$ has exactly one infinite $*$-connected component, which we denote by $\Ext(K)$; we further let
\[
\overline{K} := \Z^2\setminus \Ext(K)\,,\qquad
\Int(K) := \overline K\setminus K\,,
\]
so that $\overline K$, dubbed the \emph{$*$-hull} of $K$, is the union of $K$ and the bounded $*$-connected components of $\Z^2\setminus K$, referred to as \emph{$*$-holes} of $K$; the \emph{interior} of $K$, dubbed $\Int(K)$, is the union of these $*$-holes. 
Note that, as each $*$-hole of $K$ is a separate $*$-connected component, distinct from $\Ext(K)$, one has
\begin{equation}\label{eq:int-ext-dist}
\dist_*(\Int(K),\Ext(K))\geq 2\,.
\end{equation}
\begin{fact}\label{fact:simply}
If $V\subset\Z^2$ is finite and simply connected (i.e., $V=\overline V$), then $\overline K\subset V$ for every $K\subset V$.
\end{fact}
\begin{proof}
Let $K\subset V$ and suppose there exists $z\in\overline K \setminus V$. Since $\overline K\setminus K = \Int(K)$, we must have $z\in H$ for some $*$-hole $H$ of $K$. As $V\supset K$, the $*$-connected component of $z$ in $\Z^2\setminus V$ is contained in that of $z$ in $\Z^2\setminus K$---namely, $H$---thus it is bounded, contradicting the simple connectivity of $V$.
\end{proof}

\begin{definition}[Nesting]\label{def:nested}
Let $K_1,K_2\subset\Z^2$ be finite and nonempty. We say that $K_1$ is \emph{nested} in $K_2$, denoted $K_1\prec K_2$, if $\overline{K_1}\subset\Int(K_2)$. We write $K_1\preceq K_2$ if $K_1\prec K_2$ or $K_1=K_2$.
\end{definition}
Note that nesting is a strict partial order: it is transitive, as $\overline{K_1}\subset\Int(K_2)\subset\overline{K_2}\subset\Int(K_3)$ whenever $K_1\prec K_2\prec K_3$, and it is irreflexive, as $\Int(K)=\overline K\setminus K$. 

The following flavor of a \emph{boundary} of a finite set $K\subset\Z^2$ will be central to our construction:
\begin{equation}\label{eq:boundary-def}
\partialhat K := \big\{ y\in K \,:\, \dist_*(y,\Ext(K))= 1\big\}\,.
\end{equation}
Observe that, unlike the internal vertex boundary $\partialvtxin K$ defined above, the set $\partialhat K$ does not consider $*$-holes to be exterior to $K$, and further uses $*$-adjacency (rather than adjacency) to the exterior.
These two features, together with the feature that the exterior is the infinite $*$-connected (rather than connected) component of $\Z^2\setminus K$, crucially guarantee (\cref{lem:one-boundary,lem:two-boundaries} below) that the boundary of a connected set is connected, and that if two such boundaries are disjoint then either their $*$-hulls are disjoint or one of them is nested in the other. The examples in \cref{fig:boundary-def,fig:boundary-true} demonstrate that none of the three features can be dropped.

\begin{figure}
\centering
\begin{tikzpicture}[scale=0.3, site/.style={circle,inner sep=0pt,minimum size=4.4pt}, gsite/.style={site,draw=green!50!black,fill=green!20,thin}, bsite/.style={site,fill=blue,draw=blue}]
\draw[gray!30,very thin] (-1.5,-1.5) grid (7.5,7.5);
\foreach \x/\y in {1/1,1/2,1/3,1/4,1/5,2/1,2/5,3/1,3/5,4/1,4/5,5/1,5/2,5/3,5/4,5/5} \node[gsite] at (\x,\y) {};
\foreach \x/\y in {0/0,0/1,0/2,0/3,0/4,0/5,0/6,1/0,1/6,2/0,2/2,2/3,2/4,2/6,3/0,3/2,3/4,3/6,4/0,4/2,4/3,4/4,4/6,5/0,5/6,6/0,6/1,6/2,6/3,6/4,6/5,6/6} \node[bsite] at (\x,\y) {};
\node at (3,-3) {(a)};
\end{tikzpicture}\hspace{4em}
\begin{tikzpicture}[scale=0.3, site/.style={circle,inner sep=0pt,minimum size=4.4pt}, gsite/.style={site,draw=green!50!black,fill=green!20,thin}, bsite/.style={site,fill=blue,draw=blue}]
\draw[gray!30,very thin] (-1.5,-1.5) grid (7.5,7.5);
\foreach \x/\y in {1/1,1/5,2/1,2/2,2/3,2/4,2/5,3/1,3/2,3/3,3/4,3/5,4/1,4/2,4/3,4/4,4/5,5/1,5/2,5/3,5/4,5/5} \node[gsite] at (\x,\y) {};
\foreach \x/\y in {0/0,0/1,0/2,0/4,0/5,0/6,1/0,1/6,2/0,2/6,3/0,3/6,4/0,4/6,5/0,5/6,6/0,6/1,6/2,6/3,6/4,6/5,6/6} \node[bsite] at (\x,\y) {};
\node[bsite] at (1,3) {};
\node at (3,-3) {(b)};
\end{tikzpicture}\hspace{4em}
\begin{tikzpicture}[scale=0.3, site/.style={circle,inner sep=0pt,minimum size=4.4pt}, gsite/.style={site,draw=green!50!black,fill=green!20,thin}, b1/.style={site,fill=blue,draw=blue}, b2/.style={site,fill=violet,draw=violet}, osite/.style={site,draw=orange!70!black,fill=orange!30,thin}]
\draw[gray!30,very thin] (-4.5,-4.5) grid (4.5,4.5);
\draw[blue,thin] (-4.5,-4.5)--(4.5,4.5) (-4.5,-2.5)--(2.5,4.5);
\draw[violet,thin] (-3.5,4.5)--(4.5,-3.5) (-1.5,4.5)--(4.5,-1.5);
\foreach \x/\y in {-4/-3,-3/-2,-2/-1,-1/0,2/3,3/4} \node[gsite] at (\x,\y) {};
\foreach \x/\y in {-2/4,-1/3,2/0,3/-1,4/-2} \node[osite] at (\x,\y) {};
\foreach \x/\y in {-4/-4,-4/-2,-3/-3,-3/-1,-2/-2,-2/0,-1/-1,-1/1,0/0,1/3,2/2,2/4,3/3,4/4} \node[b1] at (\x,\y) {};
\foreach \x/\y in {-3/4,-2/3,-1/2,-1/4,0/3,1/0,2/-1,2/1,3/-2,3/0,4/-3,4/-1} \node[b2] at (\x,\y) {};
\fill[green!20] (0,1) ++(0,7.33pt) arc[start angle=90,end angle=270,radius=7.33pt] -- cycle; \fill[violet] (0,1) ++(0,-7.33pt) arc[start angle=270,end angle=450,radius=7.33pt] -- cycle; \draw[black,thin] (0,1) circle[radius=7.33pt];
\fill[blue] (0,2) ++(0,7.33pt) arc[start angle=90,end angle=270,radius=7.33pt] -- cycle; \fill[orange!30] (0,2) ++(0,-7.33pt) arc[start angle=270,end angle=450,radius=7.33pt] -- cycle; \draw[black,thin] (0,2) circle[radius=7.33pt];
\fill[blue] (1,1) ++(0,7.33pt) arc[start angle=90,end angle=270,radius=7.33pt] -- cycle; \fill[orange!30] (1,1) ++(0,-7.33pt) arc[start angle=270,end angle=450,radius=7.33pt] -- cycle; \draw[black,thin] (1,1) circle[radius=7.33pt];
\fill[green!20] (1,2) ++(0,7.33pt) arc[start angle=90,end angle=270,radius=7.33pt] -- cycle; \fill[violet] (1,2) ++(0,-7.33pt) arc[start angle=270,end angle=450,radius=7.33pt] -- cycle; \draw[black,thin] (1,2) circle[radius=7.33pt];
\node at (0,-6) {(c)};
\end{tikzpicture}
\vspace{-0.5em}
\caption{Each feature of \cref{eq:boundary-def} is needed for \cref{lem:one-boundary,lem:two-boundaries}. Sites of $K$ are circles, in blue when they are in the boundary $\partialhat K$ under its alternate definition and in green otherwise. In~(c), a second set $K'$ is in purple when sites are in its boundary and in orange otherwise. 
(a)~A~disconnected boundary when  $*$-holes are part of the exterior. (b)~A disconnected boundary when the exterior is the infinite connected (not $*$-connected) component of $\Z^2\setminus K$. (c)~Taking sites adjacent (not $*$-adjacent) to $\Ext(K)$ makes the two bands $K,K'$ have disjoint boundaries, although their $*$-hulls intersect and neither one contains the other.}
\label{fig:boundary-def}
\end{figure}
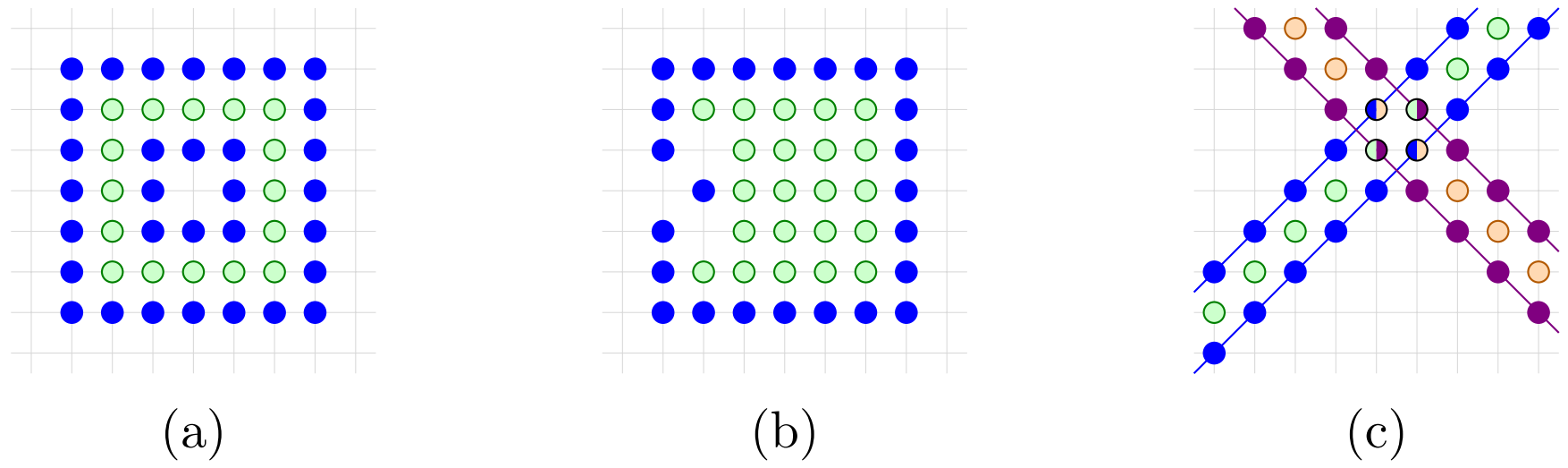

\begin{figure}
\centering
\begin{tikzpicture}[scale=0.3, site/.style={circle,inner sep=0pt,minimum size=4.4pt}, gsite/.style={site,draw=green!50!black,fill=green!20,thin}, bsite/.style={site,fill=blue,draw=blue}]
\draw[gray!30,very thin] (-1.5,-1.5) grid (7.5,7.5);
\foreach \x/\y in {1/1,1/2,1/3,1/4,1/5,2/1,2/2,2/3,2/4,2/5,3/1,3/2,3/4,3/5,4/1,4/2,4/3,4/4,4/5,5/1,5/2,5/3,5/4,5/5} \node[gsite] at (\x,\y) {};
\foreach \x/\y in {0/0,0/1,0/2,0/3,0/4,0/5,0/6,1/0,1/6,2/0,2/6,3/0,3/6,4/0,4/6,5/0,5/6,6/0,6/1,6/2,6/3,6/4,6/5,6/6} \node[bsite] at (\x,\y) {};
\node at (3,-3) {(a)};
\end{tikzpicture}\hspace{4em}
\begin{tikzpicture}[scale=0.3, site/.style={circle,inner sep=0pt,minimum size=4.4pt}, gsite/.style={site,draw=green!50!black,fill=green!20,thin}, bsite/.style={site,fill=blue,draw=blue}]
\draw[gray!30,very thin] (-1.5,-1.5) grid (7.5,7.5);
\foreach \x/\y in {3/1,3/2,3/3,3/4,3/5,4/1,4/2,4/3,4/4,4/5,5/1,5/2,5/3,5/4,5/5} \node[gsite] at (\x,\y) {};
\foreach \x/\y in {0/0,0/1,0/2,0/4,0/5,0/6,1/0,1/1,1/3,1/5,1/6,2/0,2/1,2/2,2/3,2/4,2/5,2/6,3/0,3/6,4/0,4/6,5/0,5/6,6/0,6/1,6/2,6/3,6/4,6/5,6/6} \node[bsite] at (\x,\y) {};
\node at (3,-3) {(b)};
\end{tikzpicture}\hspace{4em}
\begin{tikzpicture}[scale=0.3, site/.style={circle,inner sep=0pt,minimum size=4.4pt}, gsite/.style={site,draw=green!50!black,fill=green!20,thin}, b1/.style={site,fill=blue,draw=blue}, b2/.style={site,fill=violet,draw=violet}, osite/.style={site,draw=orange!70!black,fill=orange!30,thin}]
\draw[gray!30,very thin] (-4.5,-4.5) grid (4.5,4.5);
\foreach \x/\y in {-4/-4,-4/-3,-4/-2,-3/-3,-3/-2,-3/-1,-2/-2,-2/-1,-2/0,-1/-1,-1/0,-1/1,0/0,1/3,2/2,2/3,2/4,3/3,3/4,4/4} \node[b1] at (\x,\y) {};
\foreach \x/\y in {-3/4,-2/3,-2/4,-1/2,-1/3,-1/4,0/3,1/0,2/-1,2/0,2/1,3/-2,3/-1,3/0,4/-3,4/-2,4/-1} \node[b2] at (\x,\y) {};
\fill[blue] (0,1) ++(0,7.33pt) arc[start angle=90,end angle=270,radius=7.33pt] -- cycle; \fill[violet] (0,1) ++(0,-7.33pt) arc[start angle=270,end angle=450,radius=7.33pt] -- cycle; \draw[black,thin] (0,1) circle[radius=7.33pt];
\fill[blue] (0,2) ++(0,7.33pt) arc[start angle=90,end angle=270,radius=7.33pt] -- cycle; \fill[violet] (0,2) ++(0,-7.33pt) arc[start angle=270,end angle=450,radius=7.33pt] -- cycle; \draw[black,thin] (0,2) circle[radius=7.33pt];
\fill[blue] (1,1) ++(0,7.33pt) arc[start angle=90,end angle=270,radius=7.33pt] -- cycle; \fill[violet] (1,1) ++(0,-7.33pt) arc[start angle=270,end angle=450,radius=7.33pt] -- cycle; \draw[black,thin] (1,1) circle[radius=7.33pt];
\fill[blue] (1,2) ++(0,7.33pt) arc[start angle=90,end angle=270,radius=7.33pt] -- cycle; \fill[violet] (1,2) ++(0,-7.33pt) arc[start angle=270,end angle=450,radius=7.33pt] -- cycle; \draw[black,thin] (1,2) circle[radius=7.33pt];
\node at (0,-6) {(c)};
\end{tikzpicture}
\vspace{-0.5em}
\caption{The boundaries, per \cref{eq:boundary-def}, of the sets shown in \cref{fig:boundary-def}, with the same colors. In (a) and (b) the boundary is connected, and in (c) every site of each band is in its boundary, so the boundaries of $K$ and $K'$ intersect.}
\label{fig:boundary-true}
\end{figure}
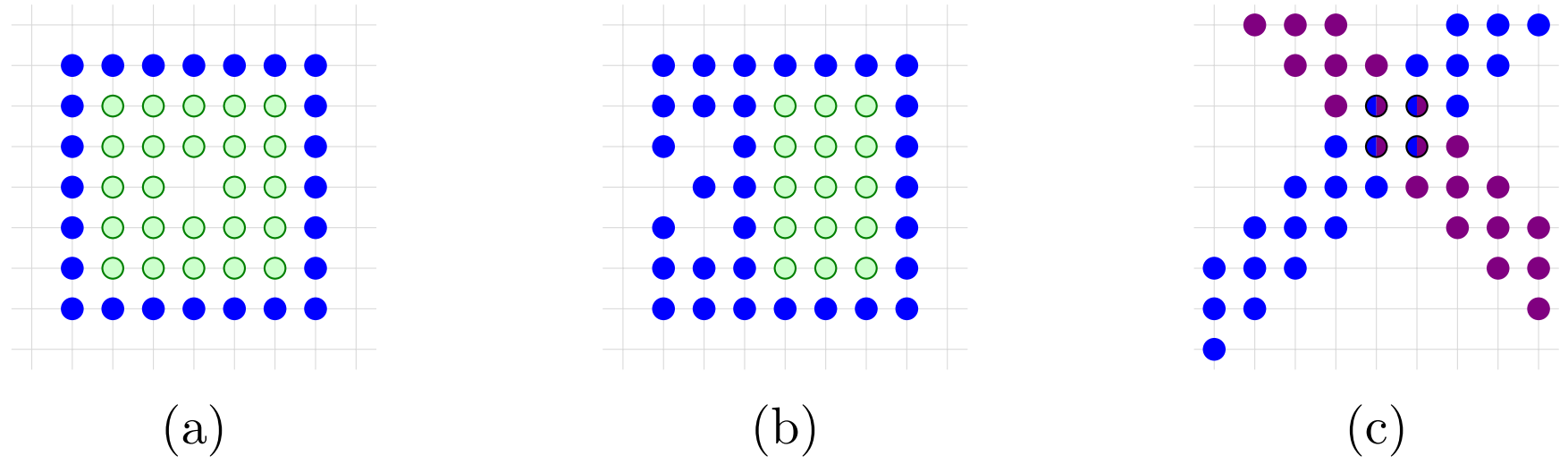

\begin{lemma}[one boundary]\label{lem:one-boundary}
Let $K\subset\Z^2$ be finite, nonempty and connected, and let $\gamma=\partialhat K$. Then the following hold:
\begin{enumerate}[(i)]
\item\label{it:one-ext} the sets $\gamma,K$ have the same exterior: $\Ext(\gamma)=\Ext(K)$ (thus $\overline\gamma=\overline K$ and $\Int(\gamma)=\overline K\setminus\partialhat K$);
\item\label{it:one-hull} the $*$-hull $\overline\gamma$ is connected and satisfies $\partialhat\overline\gamma=\gamma$; 
\item\label{it:one-conn} the set $\gamma$ is connected;
\item\label{it:one-simply} the set $\Int(\gamma)$ is simply connected.
\end{enumerate}
\end{lemma}

To prove \cref{lem:one-boundary}\labelcref{it:one-conn}, we will use the following result of Hammond (see also \cite[Thm.~4]{Timar13}), the analog in a box of a classical result of Kesten \cite[Lem.~2.23]{Kesten86} on the connectedness of boundaries.

\begin{lemma}[{\cite[Lem.~2.1]{Hammond06}}]\label{lem:hammond}
Let $Q=\llbracket1,n\rrbracket^d$ in $\Z^d$, let $\sC \subset Q$ be $*$-connected and let $x_0\in Q\setminus\sC$. Then $Y=\{y\in Q\setminus \sC:\, \mbox{$y$ is connected to $x_0$ by a path in $Q\setminus \sC$ and $\exists z\in \sC$: $ y\starsim z$}
\}$ is connected.
\end{lemma}

\begin{proof}[Proof of \cref{lem:one-boundary}]
For any finite $A\subset B$, any $*$-connected component of $B^c$ is contained in a $*$-connected component of $A^c$, so $\Ext(B)\subset\Ext(A)$; in particular, $\gamma\subset K$ implies $\Ext(K) \subset \Ext(\gamma)$. Suppose now that $X := \Ext(\gamma)\setminus\Ext(K)\neq\emptyset$, and take $u\in X$ such that 
$u\starsim v$ for some $v\in \Ext(K)$ (such a $u$ exists since $\Ext(\gamma)$ is $*$-connected). Since $u\notin \Int(K)$ by \cref{eq:int-ext-dist}, we have $u\in K$, and together with $\dist_*(u,\Ext(K))=1$, the criterion in \labelcref{eq:boundary-def} implies that $u\in \partialhat K = \gamma$, a contradiction since $\gamma\cap \Ext(\gamma)=\emptyset$. We thus arrive at \labelcref{it:one-ext}.

Having shown that $\overline\gamma=\overline K$, we now argue that it is connected. To see this, recall that $K$ is connected by assumption, let $H$ be a $*$-hole of $K$, and let $\cC$ be a connected component of $H$. If $x=(x_1,x_2)\in \cC$ is such that $x_1$ is maximal over (the finite) $\cC$, then $y=(x_1+1,x_2)$ must be in~$K$: indeed, $\cC$ is a connected component of $\Z^2\setminus K$, and $\dist(y,\cC)=1$. We have thus shown that $\cC$ is adjacent to $K$, and it follows that $\overline K$ is connected.
The fact $\partialhat\overline K=\partialhat K$ is immediate since the $*$-holes of $K$ are precluded from being $*$-adjacent to $\Ext(K)$ by \cref{eq:int-ext-dist}. This establishes \labelcref{it:one-hull}.

We turn to \labelcref{it:one-conn}, which we will infer from \cref{lem:hammond}. Let $Q$ be a box in $\Z^2$ containing every site $y$ with $\dist_*(y,\overline K)\leq1$, and let $\sC:=\Ext(K)\cap Q$. We claim that $\sC$ is $*$-connected: indeed, by the definition of $Q$, the set $\partialvtxin Q$ is $*$-connected and has $\dist_*(\partialvtxin Q,\overline K)\geq 1$; in particular, $\partialvtxin Q\subset \sC$, and so the $*$-connectedness of any two $u,v\in\sC$ is inherited from that in $\Ext(K)$ (either in $Q\setminus \partialvtxin Q$ or through the $*$-connected $\partialvtxin Q$).
Take an arbitrary $x_0\in K$ (noting $x_0\in Q\setminus \sC$). By \cref{lem:hammond}, applied to (a translate of) $\sC \subset Q$ and $x_0$, the corresponding set $Y$ from that lemma is connected. Since $\overline K\subset Q$, we have $Q\setminus\sC=Q\setminus\Ext(K)=\overline K$, which is connected (by \labelcref{it:one-hull}) and contains $x_0$. In particular, every $y\in Q\setminus \sC$ is connected to $x_0$ by a path in $Q\setminus\sC$, and so
\begin{align*} 
  Y &= \left\{ y\in \overline K\,:\;\dist_*(y,\sC) = 1\right\} = \left\{ y\in \overline K\,:\;\dist_*(y,\Ext(K)) = 1\right\}\,, 
  \end{align*}
with the second equality due to the fact $Q\supset \{y:\dist_*(y,\overline K)\leq 1\}$. Of course, $\Ext(\overline K)=\Ext(K)$ (by the definition of $\overline K$), and so, by the definition of the boundary in \labelcref{eq:boundary-def},  
\[ Y = \partialhat\overline K = \partialhat K = \gamma\]
(inferring the middle identity from \labelcref{it:one-hull}), which yields the connectedness of $\gamma$ as stated in \labelcref{it:one-conn}.

The final item \labelcref{it:one-simply} readily follows: $\Z^2\setminus\Int(\gamma)=\gamma\cup\Ext(\gamma)$ is $*$-connected (i.e., $\Int(\gamma)$ is simply connected) since $\gamma$ is connected (by \labelcref{it:one-conn}) and $*$-adjacent to $\Ext(\gamma)$; this completes the proof.
\end{proof}

The following lemma will address the nesting of two boundaries $\gamma_1 = \partialhat K_1$ and $\gamma_2=\partialhat K_2$. NB.~this is different from the nesting of the regions $K_i$ (for example, if $K_1=\{x\}$ and $K_2=B(x,1)$, then both $\overline{\gamma_1}$ and $\Int(\gamma_2)$ are the singleton $\{x\}$, whence $\gamma_1 \prec \gamma_2$, and yet $\Int(K_2)=\emptyset$ so $K_1 \nprec K_2$). 
\begin{lemma}[two boundaries]\label{lem:two-boundaries} 
  Let $K_1,K_2\subset\Z^2$ be finite, nonempty and connected, let $\gamma_1=\partialhat K_1$ and  $\gamma_2=\partialhat K_2$. With $\prec$ as per \cref{def:nested}, at most one of the following three alternatives holds:
\[
\overline{\gamma_1}\cap\overline{\gamma_2}=\emptyset\,,\qquad \gamma_1\prec\gamma_2\,,\qquad \gamma_2\prec\gamma_1\,.
\]
Moreover, if $\gamma_1\cap\gamma_2=\emptyset$ then exactly one of the three holds. Finally, if $\overline{\gamma_1}\cap\overline{\gamma_2}=\emptyset$ and $\dist(\overline{\gamma_1},\overline{\gamma_2})\leq1$ then $\dist(\gamma_1,\gamma_2)\leq1$.
\end{lemma}

\begin{proof}
For the first assertion, recall that $\gamma_i\neq\emptyset$ for each $i$, since $K_i$ is finite and nonempty. Having $\gamma_1\prec\gamma_2$ would imply that the nonempty set $\overline{\gamma_1}$ is contained in $\Int(\gamma_2)\subset\overline{\gamma_2}$, hence $\overline{\gamma_1},\overline{\gamma_2}$ are not disjoint; the same applies to $\gamma_2\prec\gamma_1$. Finally, $\gamma_1\prec\gamma_2$ and $\gamma_2\prec\gamma_1$ cannot both hold, since nesting is a strict partial order.

Suppose now that $\gamma_1\cap\gamma_2=\emptyset$ and yet none of the above three alternatives holds. We will deduce from this assumption the existence of two sites $y$ and $z$ as follows:

\begin{enumerate}[(1)]

\item\label{it:two-y} $y\in \Int(\gamma_1)\cap\gamma_2$: It suffices to find some $y\in\overline{\gamma_1}\cap \gamma_2$, as the assumption
 $\gamma_1\cap\gamma_2=\emptyset$ would imply that $y\in \overline{\gamma_1}\setminus\gamma_1=\Int(\gamma_1)$. To do so, note that, since $\gamma_1\nprec\gamma_2$, i.e.,
$\overline{\gamma_1}\not\subset\Int(\gamma_2)=\overline{\gamma_2}\setminus\gamma_2$, some site in $\overline{\gamma_1}$ must belong either to $\gamma_2$ or to $\Z^2\setminus\overline{\gamma_2}=\Ext(\gamma_2)$. In the former case we are done (having found said $y$), so assume the latter. The set $\overline{\gamma_1}$ is connected (by \cref{lem:one-boundary}\labelcref{it:one-hull}) and intersects both $\Ext(\gamma_2)$ and $\overline{\gamma_2}$ (as $\overline{\gamma_1}\cap\overline{\gamma_2}\neq\emptyset$); hence, there exists a pair $x,y\in\overline{\gamma_1}$ such that $x\in \Ext(\gamma_2)$, $y\in\overline{\gamma_2}$ and $x\sim y$. 
Observe that $\dist_*(y,\Ext(\gamma_2))=1$, and as $\Ext(\gamma_2)=\Ext(\overline{\gamma_2})$, this means that $y\in\partialhat\overline{\gamma_2}$, which is $\gamma_2$ by \cref{lem:one-boundary}\labelcref{it:one-hull}.

\item\label{it:two-z} $z\in\Ext(\gamma_1)\cap\gamma_2$: 
First, let $u\in\gamma_1\cap\Int(\gamma_2)$, the symmetric analogue of the site in \cref{it:two-y}.
Being in $\gamma_1=\partialhat K_1$, it is $*$-adjacent to some $v\in\Ext(K_1)$, which is $\Ext(\gamma_1)$ by \cref{lem:one-boundary}\labelcref{it:one-ext}. Since $u\in\Int(\gamma_2)$ has $\dist_*(u,\Ext(\gamma_2))\geq2$ by \cref{eq:int-ext-dist}, we get $v\notin\Ext(\gamma_2)$, which combines into $v\in\Ext(\gamma_1)\cap\overline{\gamma_2}$.
In particular, $\Ext(\gamma_1)\cap\overline{\gamma_2}$ is nonempty and finite.
On the other hand, $\Ext(\gamma_1)$ is $*$-connected and infinite, so there exist $z\in\Ext(\gamma_1)\cap\overline{\gamma_2}$ and $z'\in\Ext(\gamma_1)\setminus\overline{\gamma_2}$ with $z \starsim z'$. As $z\in\overline{\gamma_2}$ and $z'\in\Z^2\setminus\overline{\gamma_2}=\Ext(\overline{\gamma_2})$, we get $z\in\partialhat\overline{\gamma_2}=\gamma_2$ by \cref{lem:one-boundary}\labelcref{it:one-hull}.
\end{enumerate}
Consider the sites $y,z$ from \cref{it:two-y,it:two-z}. Since $\gamma_2$ is connected by \cref{lem:one-boundary}\labelcref{it:one-conn}, there is a path $y=w_0\sim w_1\sim\cdots\sim w_m=z$ in $\gamma_2$. Noting that $z\notin \Int(\gamma_1)$, let $w_j$ be the first site (minimal index $j\geq 1$) along this path that does not belong to $\Int(\gamma_1)$. Being adjacent to $w_{j-1}\in\Int(\gamma_1)$, the site $w_j$ cannot be in $\Ext(\gamma_1)$ by \cref{eq:int-ext-dist}, so it must be that $w_j\in\gamma_1$, contradicting the assumption that $\gamma_1\cap\gamma_2=\emptyset$. We therefore conclude that one, and exactly one, of the three alternatives holds.

For the last assertion, we have $\overline{\gamma_1}\cap\overline{\gamma_2}=\emptyset$ and, as $\dist(\overline{\gamma_1},\overline{\gamma_2})\leq1$, may take $u_1\in\overline{\gamma_1}$ and $u_2\in\overline{\gamma_2}$ such that $u_1\sim u_2$. Then $u_2\notin\overline{\gamma_1}$, so $u_2\in\Ext(\gamma_1)$, and $u_1\in\overline{\gamma_1}$ has $\dist_*(u_1,\Ext(\gamma_1))=1$; thus, $u_1\in\partialhat\overline{\gamma_1}=\gamma_1$ by \cref{lem:one-boundary}\labelcref{it:one-hull}. Similarly, $u_2\in\gamma_2$, so $\dist(\gamma_1,\gamma_2)\leq1$, as claimed.
\end{proof}

\subsection{Proof of \texorpdfstring{\cref{thm:main}}{the main theorem}}\label{subsec:proof}
We are now ready to define the polymer model.
\begin{definition}[Polymers]
With \cref{eq:boundary-def} in mind, define the set of \emph{polymer supports} in $\Lambda$ to be
\[
\Gamma := \big\{ \partialhat K \,:\, K\subset\Lambda\mbox{ is nonempty and connected}\big\}\,.
\]
Two elements $\gamma_1,\gamma_2\in\Gamma$ are \emph{compatible}, denoted $\gamma_1\perp\gamma_2$, if
\begin{equation}\label{eq:compatible-def}
\gamma_1\prec\gamma_2\quad\mbox{or}\quad \gamma_2\prec\gamma_1\quad\mbox{or}\quad \dist(\gamma_1,\gamma_2)\geq2
\end{equation}
(either one is nested in the other, as per \cref{def:nested}, or they neither intersect nor are adjacent). 
A~\emph{polymer} is a pair $\pi=(\gamma,\eta)$ for $\gamma\in\Gamma$ (its support) and $\eta\in \scS^\gamma$ (spin values on the support); let~$\gamma_\pi$ and $\eta_\pi$ denote its components, and let $\cP:=\{(\gamma,\eta):\gamma\in\Gamma,\ \eta\in\scS^\gamma\}$ be the set of polymers. The $*$-hull of a polymer is that of its support; for $V\subset\Lambda$, let $\cP_V\subset\cP$ denote the set of polymers whose $*$-hulls are contained in $V$.
Two polymers are compatible if their supports are; a polymer is nested in another if its support is nested in that of the other; and a family $\cF$ of polymers is compatible if its elements are pairwise compatible.
\end{definition}
Note that $\gamma\not\perp\gamma$ for every $\gamma\in\Gamma$ (as $\gamma$ is not nested in itself and is at distance $0$ from itself), so a compatible family contains at most one polymer with a given support. Both $\Gamma$ and the compatibility relation are purely geometric notions: they depend on $\Lambda$, but neither on $\fg$ nor on the potentials~$\scV$.
\begin{remark}\label{rem:polymers}
The set $\Gamma$, relation $\perp$ and set $\cP$ are those of \cref{thm:main}, and satisfy the properties given there: every $\gamma\in\Gamma$ is a connected subset of $\Lambda$ by \cref{lem:one-boundary}\labelcref{it:one-conn}; if $\gamma_1\perp\gamma_2$ then $\gamma_1\cap\gamma_2=\emptyset$, since $\gamma_1\prec\gamma_2$ implies $\gamma_1\subset\overline{\gamma_1}\subset\overline{\gamma_2}\setminus\gamma_2$; and if $\gamma_1\not\perp\gamma_2$ then $\dist(\gamma_1,\gamma_2)\leq1$ by \cref{eq:compatible-def}.
\end{remark}

\begin{definition}[Maximal polymer, antichain]\label{def:max-family}
For a family of polymers $\cF$, we say $\pi\in\cF$ is \emph{maximal} if it is not nested in any other $\pi'\in\cF$. Call $\cF$ an \emph{antichain} if all its elements are maximal.
\end{definition}
\begin{observation}\label{obs:nesting}
The $*$-hulls of all maximal elements of a compatible family of polymers $\cF$ have pairwise distances at least two; every other element of $\cF$ is nested in exactly one maximal element. In particular, a family $\cF\subset\cP$ is a compatible antichain if and only if the $*$-hulls of its elements are pairwise at distance at least two.
\end{observation}
\begin{proof}
Let $\cF = \{\pi_i\}$, and within this proof, write $\gamma_i = \gamma_{\pi_i}$ for brevity. We first argue that if $\pi_1,\pi_2\in\cF$ are distinct and neither one is nested in the other, then $\dist(\overline{\gamma_1},\overline{\gamma_2})\geq2$. Indeed, as $\gamma_1\perp\gamma_2$, in this case \cref{eq:compatible-def} reads $\dist(\gamma_1,\gamma_2)\geq2$. In particular, $\gamma_1\cap\gamma_2=\emptyset$, so exactly one of the three alternatives in \cref{lem:two-boundaries} holds; since, again, $\gamma_1\nprec\gamma_2$ and $\gamma_2\nprec\gamma_1$, it must be that $\overline{\gamma_1}\cap\overline{\gamma_2}=\emptyset$.
Having $\dist(\overline{\gamma_1},\overline{\gamma_2})\leq1$ would then imply that $\dist(\gamma_1,\gamma_2)\leq1$ by the last assertion of that lemma, a contradiction. Applying this to two distinct maximal elements of $\cF$ (neither of which is nested in the other, by definition) establishes the first assertion.

For the second assertion, let $\pi_0\in\cF$ be an element that is not maximal. Since nesting is a strict partial order and the family $\cF$ is finite (as $\Lambda$ is), the element $\pi_0$ is nested in some maximal element of $\cF$. Moreover, if $\pi_0$ were nested in two distinct maximal elements $\pi_1,\pi_2$, then the nonempty $\overline{\gamma_0}$ would be contained in $\overline{\gamma_1}\cap\overline{\gamma_2}$, which is empty by the first assertion.

For the last assertion, its `only if' direction is the first assertion, and the converse is immediate: if the $*$-hulls of the elements of $\cF$ are pairwise at distance at least two, then so are their supports, so $\cF$ is compatible by \cref{eq:compatible-def}, and no element is nested in another, as $\gamma_1\prec\gamma_2$ would give $\overline{\gamma_1}\subset\overline{\gamma_2}$. This concludes the proof.
\end{proof}

For $(\gamma,\eta)\in\cP$, let $\mu_{\Int(\gamma)}^\eta$ be the Gibbs measure on $\Int(\gamma)$ with boundary condition $\eta$ (that is, $\mu_{\Int(\gamma)}^\tau$ for any boundary condition $\tau\in \scS^{\Z^2\setminus\Int(\gamma)}$ extending $\eta$); it is well defined since $\partialvtxout\Int(\gamma)\subset\gamma$ by \cref{eq:int-ext-dist} (recall that $\mu_V^\tau$ depends on $\tau$ only through $\tau\restriction_{\partialvtxout V}$), and it is defined whenever $Z_{\Int(\gamma)}^\eta>0$ (which is the case whenever $\mu_V^\tau(\sigma\restriction_\gamma=\eta)>0$ for some $V\supset\overline\gamma$ and $\tau\in\scS^{\Z^2\setminus V}$; recall that we allow hard constraints, without which these conditions would be automatic).

Let $(\gamma_i,\eta_i)$ be polymers whose $*$-hulls are pairwise disjoint and contained in $V\subset\Lambda$. Then
\[
\dist_*(\Int(\gamma_i),\Int(\gamma_j))\geq2\qquad\mbox{for all $i\neq j$}\,,
\]
since $\Int(\gamma_i)\subset\overline{\gamma_i}\subset\Z^2\setminus\overline{\gamma_j}=\Ext(\gamma_j)$, whereas we have $\dist_*(\Int(\gamma_j),\Ext(\gamma_j))\geq2$ by \cref{eq:int-ext-dist}.
Applying \cref{lem:markov} with $U=\bigcup_i\Int(\gamma_i)$ for these polymers, and summing over $\xi$ with $\xi\restriction_{\gamma_i}=\eta_i$, we obtain that, for all functions $F_i$ on $\scS^{\Int(\gamma_i)}$,
\begin{equation}\label{eq:markov-factorization}
\E_V^\tau\bigg[\prod_i \one_{\{\sigma\restriction_{\gamma_i}=\eta_i\}}F_i(\sigma\restriction_{\Int(\gamma_i)})\bigg] = \mu_V^\tau\big(\forall i,\, \sigma\restriction_{\gamma_i}=\eta_i\big)\prod_i\E_{\Int(\gamma_i)}^{\eta_i}[F_i]\,,
\end{equation}
where we interpret the right-hand side as $0$ if the probability that appears there vanishes (the measures $\mu_{\Int(\gamma_i)}^{\eta_i}$ need not be defined), while if it is positive then $Z_{\Int(\gamma_i)}^{\eta_i}>0$ for all $i$ by \cref{lem:markov}.

\begin{definition}[Patterns]\label{def:pattern}
A \emph{pattern} is a triple $\fp=(x,r,\omega)$ with $x\in\Lambda$ and $r\geq0$ such that $B(x,r)\subset\Lambda$, and $\omega\in \scS^{B(x,r)}$; we write $B(\fp):=B(x,r)$, $r(\fp):=r$, and set
\begin{equation}\label{eq:f-def}
f(\fp) := e^{g_{x,r}(\omega)}-1\,,\qquad
\one_\fp(\sigma) := \one_{\{\sigma\restriction_{B(x,r)}=\omega\}}\,.
\end{equation}
A set $\cC$ of patterns is \emph{compatible} if any two of its elements agree on the intersection of their balls. In that case we let $B(\cC):=\bigcup_{\fp\in\cC}B(\fp)$ and let $\omega_\cC\in \scS^{B(\cC)}$ be the common configuration, so that
\begin{equation}\label{eq:indicator-product}
\one_{\{\sigma\restriction_{B(\cC)}=\omega_\cC\}} = \prod_{\fp\in\cC}\one_\fp(\sigma)\,,
\end{equation}
whereas for an incompatible $\cC$ the right-hand side of \cref{eq:indicator-product} vanishes identically. For $V\subset\Lambda$, let $\fC(V)$ denote the set of compatible collections $\cC$ with $B(\cC)\subset V$ (including the empty collection, in which case $B(\cC)=\emptyset$ and the left-hand side of \cref{eq:indicator-product} is $1$).
\end{definition}
Two immediate observations on the notions in \cref{def:pattern} will be useful in our analysis. First, by the bound on $g_{x,r}$ from \cref{eq:g-decay} and the facts that $|e^u-1|\leq 2|u|$ for $|u|\leq 1$,
\begin{equation}
\label{eq:f-bound}
|f(\fp)|\leq e^{-(\Cst-1) (r(\fp)+1)}\qquad\mbox{for all $\fp$}\,.
\end{equation}
Second, a compatible collection of patterns $\cC$ contains at most one pattern with a given pair $(x,r)$. 

We prove \cref{eq:polymer-rep} in the following form, simultaneously for every simply connected $V\subset\Lambda$ and $\tau\in \scS^{\Z^2\setminus V}$ with $Z_V^\tau>0$ (for $V=\Lambda$, it is precisely \cref{eq:polymer-rep}), with a single family of weights $\sw$:
\begin{equation}\label{eq:ZD-polymer}
Z_V^{\tau,\fg} = \sum_{\substack{\cF\subset\cP_V\\ \text{compatible}}}\E_V^\tau\bigg[\prod_{(\gamma,\eta)\in\cF}\sw(\gamma,\eta)\,\one_{\{\sigma\restriction_\gamma=\eta\}}\bigg]\,,
\end{equation}
where, by the definition of $\cP_V$, the $*$-hulls of the polymers in $\cF$ are contained in~$V$ (as $V$ is simply connected, this is equivalent to their supports being contained in~$V$, by \cref{fact:simply}), and $\sw$ is defined in \cref{eq:w-def} below, depending on neither $V$ nor $\tau$, which is crucial for the proof. 

Our argument consists of four steps: the expansion \cref{eq:ZD-expansion} of $Z_V^{\tau,\fg}$ into patterns; the identity \cref{eq:peeling}, which peels off outermost boundaries of a collection of patterns; iterating the latter---relying on the above uniformity in the induction---yielding \cref{eq:ZD-polymer}; and the bound in \cref{eq:w-bound}.

\begin{enumerate}[label=\textbf{Step~\arabic*.}, ref=\arabic*, wide=0pt, itemsep=1ex, listparindent=\parindent, parsep=0pt]

\item\label[step]{st:expansion} \emph{(Expansion into patterns.)}
Fix $V\subset\Lambda$ and $(x,r)$ with $B(x,r)\subset V$. As per the definitions in \cref{eq:f-def}, we may write $f(x,r,\sigma\restriction_{B(x,r)})$ as a sum over $f(x,r,\omega)\,\one_{(x,r,\omega)}(\sigma)$, so that
\[
e^{g_{x,r}(\sigma\restriction_{B(x,r)})} = 1+\sum_{\omega\in \scS^{B(x,r)}} f(x,r,\omega)\,\one_{(x,r,\omega)}(\sigma)\,.
\]
Taking a product over all $(x,r)$ with $B(x,r)\subset V$, expanding, and using \cref{eq:indicator-product} while noting that collections with two patterns on the same ball are incompatible, we get (recalling \cref{eq:HD-def}) that
\begin{equation}\label{eq:Hg-expansion}e^{\cH_V^\fg(\sigma)} = \sum_{\cC\in\fC(V)}\one_{\{\sigma\restriction_{B(\cC)}=\omega_\cC\}}\prod_{\fp\in\cC} f(\fp)\,,\end{equation}
and hence
\begin{equation}\label{eq:ZD-expansion}
Z_V^{\tau,\fg} = \sum_{\cC\in\fC(V)}\mu_V^\tau\big(\sigma\restriction_{B(\cC)}=\omega_\cC\big)\prod_{\fp\in\cC} f(\fp)\,.
\end{equation}
(The sums in \cref{eq:Hg-expansion,eq:ZD-expansion} are finite as there are finitely many pairs $(x,r)$ with $B(x,r)\subset V$, so no convergence issue arises here or in the rearrangements below.)

\item\label[step]{st:peeling} \emph{(Outermost boundaries.)}
Two patterns $\fp,\fp'$ are \emph{linked} if $\dist(B(\fp),B(\fp'))\leq1$, i.e., their balls intersect or contain adjacent sites. Let $\cC\in\fC(V)$, and write $\cC=\cD_1\cup\cdots\cup\cD_m$ where the $\cD_i$ are the connected components of patterns in the graph where adjacency is given by this relation (in other words, the $\cD_i$ are equivalence classes of $\cC$ under the transitive closure of this relation). Let $B_i:=B(\cD_i)$, noting it is connected, being a union of balls $B(x,r)$ any two of which are joined by a chain of intersecting or adjacent balls (each ball is a box, thus connected), and that $\dist(B_i,B_j)\geq 2$ for all $i\neq j$. Assume henceforth that $V\subset\Lambda$ is simply connected (recall that $\Lambda$ is simply connected by hypothesis), whence $\overline{B_i}\subset V$ by \cref{fact:simply}. Set
\begin{equation}\label{eq:linked-components}
\gamma_i := \partialhat B_i\in\Gamma\,,\qquad \eta_i := \omega_\cC\restriction_{\gamma_i}\,,\qquad i=1,\ldots,m\,.
\end{equation}
Recall that $\overline{\gamma_i}=\overline{B_i}$ for every $i$ by \cref{lem:one-boundary}\labelcref{it:one-ext}. We say that the component $\cD_i$ is nested in $\cD_j$ if $\gamma_i\prec\gamma_j$, that is, if $\overline{B_i}\subset\Int(\gamma_j)$, and call $\cD_i$ \emph{outermost} if it is not nested in any other component. The polymers $(\gamma_i,\eta_i)$ form a compatible family, as $\gamma_i\subset B_i$ and $\gamma_j\subset B_j$ have $\dist(\gamma_i,\gamma_j)\geq 2$ for $i\neq j$ as stated above, and $\cD_i$ is outermost if and only if $(\gamma_i,\eta_i)$ is a maximal element of this family; hence, by \cref{obs:nesting}, the outermost components have $*$-hulls pairwise at distance at least two, and every component $\cD_i$ satisfies $\gamma_i\preceq\gamma_j$ for exactly one outermost component $\cD_j$. Thus,
\[
\cF_*(\cC) := \big\{(\gamma_i,\eta_i) \,:\, \cD_i\mbox{ is outermost}\big\}
\]
is a compatible antichain in $\cP_V$ by \cref{obs:nesting}, and $\cC$ is the disjoint union of the sub-collections
\[
\cC_i := \bigcup_{j\,:\,\gamma_j\preceq\gamma_i}\cD_j\qquad(\cD_i\mbox{ is outermost})\,;
\]
each $\cC_i$ is a compatible collection with $B(\cC_i)\subset\overline{B_i}=\overline{\gamma_i}$ and $\cD_i$ as its unique outermost component.

For $(\gamma,\eta)\in\cP$, define the set $\fC_{\mathfrak c}(\gamma,\eta)\subset\fC(\overline\gamma)$ of \emph{constrained} collections to consist of $\cC\in\fC(\overline\gamma)$ which have exactly one outermost component, and this component $\cD$ further satisfies $\partialhat B(\cD)=\gamma$ and $\omega_\cC\restriction_\gamma=\eta$ (i.e., the constraint is to have $\gamma$, with spins $\eta$, as the outermost boundary of $\cC$). Whenever $Z_{\Int(\gamma)}^\eta>0$, we define the constrained partition function to be
\begin{equation}\label{eq:Ec-def}
E_{\mathfrak c}(\gamma,\eta) := \sum_{\cC\in\fC_{\mathfrak c}(\gamma,\eta)} \prod_{\fp\in\cC}f(\fp)\, \mu_{\Int(\gamma)}^\eta\Big(\sigma\restriction_{B(\cC)\setminus\gamma}=\omega_\cC\restriction_{B(\cC)\setminus\gamma}\Big)\,,
\end{equation}
which is well defined since $B(\cC)\setminus\gamma\subset\Int(\gamma)$; if $Z_{\Int(\gamma)}^\eta=0$, set $E_{\mathfrak c}(\gamma,\eta):=0$. 

We claim that, for every simply connected $V\subset\Lambda$ and every $\tau$ with $Z_V^\tau>0$,
\begin{equation}\label{eq:peeling}
Z_V^{\tau,\fg} = \sum_{\substack{\cF_*=\{(\gamma_i,\eta_i)\}\subset\cP_V \\\text{compatible}\\\text{antichain}}} \mu_V^\tau\big(\forall i,\, \sigma\restriction_{\gamma_i}=\eta_i\big)\prod_{i}E_{\mathfrak c}(\gamma_i,\eta_i)\,.
\end{equation}
To prove this, we will reorganize the sum over $\cC\in\fC(V)$ in \cref{eq:ZD-expansion} via the map $\cC\mapsto(\cF_*(\cC),(\cC_i)_i)$, which the next argument establishes as a bijection from $\fC(V)$ onto the set of pairs consisting of a compatible antichain $\cF_*=\{(\gamma_i,\eta_i)\}\subset\cP_V$ and an element $\cC_i\in \fC_{\mathfrak c}(\gamma_i,\eta_i)$ for each $i$:
\begin{enumerate}[(a)]
\item As argued above, every $\cC\in\fC(V)$ is the disjoint union of the sub-collections $\cC_i\in\fC_{\mathfrak c}(\gamma_i,\eta_i)$ over the outermost components $\cD_i$, where $\cF_*(\cC)=\{(\gamma_i,\eta_i)\}$ is a compatible antichain in $\cP_V$.
\item Conversely, let $\cF_*=\{(\gamma_i,\eta_i)\}\subset\cP_V$ be a compatible antichain, and let $\cC_i\in\fC_{\mathfrak c}(\gamma_i,\eta_i)$ for each~$i$. The balls comprising $\cC_i$ are contained in $\overline{\gamma_i}$, and $\dist(\overline{\gamma_i},\overline{\gamma_j})\geq 2$ for $i\neq j$ by \cref{obs:nesting}. Hence, $\cC := \bigcup \cC_i$ is compatible, belongs to $\fC(V)$, and has no linked pair of patterns from distinct $\cC_i,\cC_j$; thus, the linked components of $\cC$ are those of the $\cC_i$. Moreover, every component of $\cC_i$ has its $*$-hull in $\overline{\gamma_i}$ by \cref{fact:simply}, as $*$-hulls are simply connected, so components from different $\cC_i$ have disjoint $*$-hulls and neither is nested in the other. We conclude that a component of $\cC$ is outermost if and only if it is outermost in its $\cC_i$, whence $\cF_*(\cC)=\cF_*$ and the sub-collections of $\cC$ are the $\cC_i$.
\end{enumerate}
Having established the bijection, we may rearrange the sum in \cref{eq:ZD-expansion}, writing each $\cC\in\fC(V)$ as the union of its sub-collections $\cC_i\in\fC_{\mathfrak c}(\gamma_i,\eta_i)$, with $(\gamma_i,\eta_i)$ the elements of $\cF_*(\cC)$, and get 
\[
Z_V^{\tau,\fg} = \sum_{\substack{\cF_*=\{(\gamma_i,\eta_i)\}\subset\cP_V \\\text{compatible}\\\text{antichain}}}\ \sum_{\substack{(\cC_i)_i\,:\\ \cC_i\in\fC_{\mathfrak c}(\gamma_i,\eta_i)}} \mu_V^\tau\big(\sigma\restriction_{B(\bigcup \cC_i)}=\omega_{\bigcup\cC_i}\big)\prod_{i}\,\prod_{\fp\in\cC_i} f(\fp)\,.
\]
Now, by \cref{eq:markov-factorization} applied with $F_i$ the indicator of $\{\sigma\restriction_{B(\cC_i)\setminus\gamma_i}=\omega_{\cC_i}\restriction_{B(\cC_i)\setminus\gamma_i}\}$, a function of $\sigma\restriction_{\Int(\gamma_i)}$,
\[
\mu_V^\tau\big(\sigma\restriction_{B(\bigcup\cC_i)}=\omega_{\bigcup\cC_i}\big) = \mu_V^\tau\big(\forall i,\, \sigma\restriction_{\gamma_i}=\eta_i\big)\prod_{i}\mu_{\Int(\gamma_i)}^{\eta_i}\Big(\sigma\restriction_{B(\cC_i)\setminus\gamma_i}=\omega_{\cC_i}\restriction_{B(\cC_i)\setminus\gamma_i}\Big)\,,
\]
where the right-hand side is interpreted as $0$ if the $\mu_V^\tau$-probability that appears there vanishes, as we did in \cref{eq:markov-factorization}. Inserting this into the previous display, the sum over the sub-collections $\cC_i$ factorizes into $\prod_i E_{\mathfrak c}(\gamma_i,\eta_i)$ by \cref{eq:Ec-def}, which establishes \cref{eq:peeling}.

\item\label[step]{st:iteration} \emph{(Iteration.)}
For $(\gamma,\eta)\in\cP$ with $Z_{\Int(\gamma)}^\eta>0$, define the unconstrained partition function of the interior of $\gamma$ with boundary condition $\eta$ to be
\begin{equation}\label{eq:Ei-def}
E_{\mathfrak i}(\gamma,\eta) := Z_{\Int(\gamma)}^{\eta,\fg} = \E_{\Int(\gamma)}^\eta\Big[e^{\cH_{\Int(\gamma)}^\fg(\sigma)}\Big] > 0\,,
\end{equation}
and define $\sw:\cP\to\R$ by
\begin{equation}\label{eq:w-def}
\sw(\gamma,\eta) := \frac{E_{\mathfrak c}(\gamma,\eta)}{E_{\mathfrak i}(\gamma,\eta)}\,,\qquad\mbox{and}\qquad \sw(\gamma,\eta):=0 \mbox{ if } Z_{\Int(\gamma)}^\eta=0\,.
\end{equation}
We prove \cref{eq:ZD-polymer} by induction on $|V|$, starting from the trivial $V=\emptyset$ where both sides equal~$1$. Let $V\neq\emptyset$ be simply connected, and apply \cref{eq:peeling}. Consider a family $\cF_*=\{(\gamma_i,\eta_i)\}$ appearing there for which $\mu_V^\tau(\forall i,\, \sigma\restriction_{\gamma_i}=\eta_i)>0$, and consider $(\gamma_i,\eta_i)\in\cF_*$ for some $i$. Then $Z_{\Int(\gamma_i)}^{\eta_i}>0$; moreover, $\Int(\gamma_i)$ is simply connected by \cref{lem:one-boundary}\labelcref{it:one-simply}, and $|\Int(\gamma_i)|<|V|$ since $\overline{\gamma_i}\subset V$ and $\gamma_i\neq\emptyset$. Thus, the induction hypothesis, applied to $\Int(\gamma_i)$ with boundary condition $\eta_i$, yields
\[ Z_{\Int(\gamma_i)}^{\eta_i,\fg} = \E_{\Int(\gamma_i)}^{\eta_i}\big[\Xi_{\Int(\gamma_i)}\big]\qquad\mbox{where}\qquad 
\Xi_{U}(\sigma):=\sum_{\substack{\cF\subset\cP_U \\\text{compatible}}}\ \prod_{(\gamma,\eta)\in\cF}\sw(\gamma,\eta)\,\one_{\{\sigma\restriction_{\gamma}=\eta\}}\,.
\]
In particular, by \cref{eq:Ei-def,eq:w-def}, we have
\[
E_{\mathfrak c}(\gamma_i,\eta_i) = \sw(\gamma_i,\eta_i)\,Z_{\Int(\gamma_i)}^{\eta_i,\fg} = \sw(\gamma_i,\eta_i)\,\E_{\Int(\gamma_i)}^{\eta_i}\big[\Xi_{\Int(\gamma_i)}\big]\,.\]
Plugging this into \cref{eq:peeling} and using the domain Markov property as per \cref{eq:markov-factorization} once more, now with $F_i:=\Xi_{\Int(\gamma_i)}$ (while noting that $\Xi_{\Int(\gamma_i)}$ is a function of $\sigma\restriction_{\Int(\gamma_i)}$), we obtain that
\[
Z_V^{\tau,\fg} = \sum_{\substack{\cF_*=\{(\gamma_i,\eta_i)\}\subset\cP_V \\\text{compatible}\\\text{antichain}}}\E_V^\tau\bigg[\prod_i \sw(\gamma_i,\eta_i)\one_{\{\sigma\restriction_{\gamma_i}=\eta_i\}}\prod_i \Xi_{\Int(\gamma_i)}(\sigma)\bigg]\,.
\]
Expanding each $\Xi_{\Int(\gamma_i)}$ according to its definition, we get
\[
Z_V^{\tau,\fg} = \sum_{\substack{\cF_*=\{(\gamma_i,\eta_i)\}\subset\cP_V \\\text{compatible}\\\text{antichain}}}\ \sum_{\substack{(\cF_i)_i\,:\\ \cF_i\subset \cP_{\Int(\gamma_i)}\\\text{compatible}}}\E_V^\tau\bigg[\prod_{(\gamma,\eta)\in\cF_*\cup\bigcup\cF_i}\sw(\gamma,\eta)\,\one_{\{\sigma\restriction_\gamma=\eta\}}\bigg]\,,
\]
The induction step will be concluded once we show that the map $(\cF_*,(\cF_i)_i)\mapsto \cF := \cF_*\cup\bigcup \cF_i$ is a bijection onto the compatible families $\cF\subset \cP_V$:
\begin{enumerate}[(a)]
\item\label{it:union-compatible} The union $\cF$ is a compatible family in $\cP_V$. To see this, first observe that $\cF_i\subset\cP_{\Int(\gamma_i)}\subset\cP_V$, so $\cF\subset\cP_V$. Moreover, the inclusion $\cF_i\subset\cP_{\Int(\gamma_i)}$ says precisely that every $\pi\in\cF_i$ is nested in $(\gamma_i,\eta_i)$. Hence, writing $\pi_i:=(\gamma_i,\eta_i)$ for the $i$-th element of $\cF_*$, we have $\cF=\bigcup_i\big(\cF_i\cup\{\pi_i\}\big)$, where, for each fixed $i$, the polymers of $\cF_i\cup\{\pi_i\}$ are pairwise compatible and have supports in $\overline{\gamma_i}$, while $\dist(\overline{\gamma_i},\overline{\gamma_j})\geq2$ for $i\neq j$ by \cref{obs:nesting}; thus $\cF$ is compatible.
\item Every compatible family $\cF\subset\cP_V$ has exactly one preimage. To see this, let $\cF_*=\{(\gamma_i,\eta_i)\}$ be the set of maximal elements of $\cF$, and let $\cF_i$ consist of the elements of $\cF$ nested in $(\gamma_i,\eta_i)$. By \cref{obs:nesting}, $\cF_*$ is a compatible antichain and the $\cF_i$ partition $\cF\setminus\cF_*$; each $\cF_i$ is compatible, and $\cF_i\subset\cP_{\Int(\gamma_i)}$ by the definition of nesting; thus, $(\cF_*,(\cF_i)_i)$ is a preimage of~$\cF$. It is the only one: for any preimage $(\cF_*',(\cF_i')_i)$ of $\cF$, every element of $\cF_i'$ is nested in the $i$-th element of $\cF_*'$, and $\cF_*'$ is an antichain, so the maximal elements of $\cF$ are exactly those of $\cF_*'$, i.e., $\cF_*'=\cF_*$, and $\cF_i'\subset\cF_i$ for every $i$, whence $\cF_i'=\cF_i$ as both families partition $\cF\setminus\cF_*$.
\end{enumerate} This completes the induction step, establishing \cref{eq:ZD-polymer} and thereby also \cref{eq:polymer-rep}.

\item\label[step]{st:bound} \emph{(The weight bound.)}
Fix $(\gamma,\eta)\in\cP$ with $Z_{\Int(\gamma)}^\eta>0$, so that $\sw(\gamma,\eta)=E_{\mathfrak c}(\gamma,\eta)/E_{\mathfrak i}(\gamma,\eta)$ by \cref{eq:w-def}. Letting
\[
\cI_\gamma := \big\{\fp \,:\; B(\fp)\cap\gamma\neq\emptyset\big\}\,,
\]
we will rewrite the sum in \cref{eq:Ec-def} via the map $\cC\mapsto(\cC_\gamma,\cC_{\mathrm{in}})$, where
\[
\cC_\gamma := \cC \cap \cI_\gamma \,,\qquad \cC_{\mathrm{in}} := \cC\setminus\cI_\gamma\,.
\]
We argue that this is a bijection from $\fC_{\mathfrak c}(\gamma,\eta)$ onto the set of pairs $(\cC_\gamma,\cC_{\mathrm{in}})\in \fC(\overline\gamma) \times \fC(\Int(\gamma))$ where $\cC_\gamma$ satisfies $\cC_\gamma\subset\cI_\gamma$, $\gamma\subset B(\cC_\gamma)$ and $\omega_{\cC_\gamma}\restriction_\gamma=\eta$, and $\cC_{\mathrm{in}}$ is compatible with $\cC_\gamma$. Since the map $X\mapsto (X\cap\cI_\gamma, X\setminus \cI_\gamma)$ is a bijection from the sets of patterns onto the pairs $(X_1,X_2)$ of such sets with $X_1\subset\cI_\gamma$ and $X_2\cap\cI_\gamma=\emptyset$ (with inverse $(X_1,X_2)\mapsto X_1\cup X_2$), and every pair as above is of this form (the balls comprising $\cC_{\mathrm{in}}$ being contained in $\Int(\gamma)$), it suffices to verify that this map sends $\fC_{\mathfrak c}(\gamma,\eta)$ into the set of pairs above, and that its inverse sends every such pair into $\fC_{\mathfrak c}(\gamma,\eta)$:
\begin{enumerate}[(a)]
\item\label{it:split-map-1} Let $\cC\in\fC_{\mathfrak c}(\gamma,\eta)$. By the definition of $\fC_{\mathfrak c}(\gamma,\eta)$, the collection of patterns $\cC$ has a unique outermost component $\cD$, with $\partialhat B(\cD)=\gamma$ and $\omega_\cC\restriction_\gamma=\eta$. Since $\gamma=\partialhat B(\cD)\subset B(\cD)$, every $x\in\gamma$ belongs to $B(\fp)$ for some $\fp\in\cD$, which also implies $\fp\in \cI_\gamma$; hence $\gamma\subset B(\cC_\gamma)\subset B(\cC)\subset\overline\gamma$, and $\omega_{\cC_\gamma}\restriction_\gamma=\omega_\cC\restriction_\gamma=\eta$. Moreover, every $\fp\in\cC_{\mathrm{in}}$ has $B(\fp)\subset B(\cC)\subset\overline\gamma$ and, as $\fp\notin\cI_\gamma$, $B(\fp)\cap\gamma=\emptyset$; thus $B(\fp)\subset\overline\gamma\setminus\gamma=\Int(\gamma)$, whence $\cC_{\mathrm{in}}\in\fC(\Int(\gamma))$. Finally, $\cC_\gamma$ and $\cC_{\mathrm{in}}$ are compatible with each other, their union $\cC$ being compatible.

\item\label{it:split-map-2} Let $(\cC_\gamma,\cC_{\mathrm{in}})$ be a pair as above, and let $\cC := \cC_\gamma\cup\cC_{\mathrm{in}}$. Note that $\cC$ is compatible with $B(\cC)\subset\overline\gamma$, so $\cC\in\fC(\overline\gamma)$, and $\omega_\cC\restriction_\gamma=\omega_{\cC_\gamma}\restriction_\gamma=\eta$. The patterns of $\cC_\gamma$ all belong to a single linked component~$\cD$ of $\cC$: each $x\in\gamma$ belongs to $B(\fp_x)$ for some $\fp_x\in\cC_\gamma$, as $\gamma\subset B(\cC_\gamma)$; the patterns $\fp_{x_1},\fp_{x_2}$ are linked or equal whenever $x_1\sim x_2$, hence all the $\fp_x$ belong to one component, $\gamma$ being connected (\cref{lem:one-boundary}\labelcref{it:one-conn}); and every $\fp\in\cC_\gamma$ has $B(\fp)\cap B(\fp_x)\neq\emptyset$ for some $x\in\gamma$, as~$\cC_\gamma\subset\cI_\gamma$. As $\cC_\gamma\subset\cD\subset\cC$, we have $\gamma\subset B(\cC_\gamma)\subset B(\cD)\subset B(\cC)\subset\overline\gamma$, so that $\Ext(B(\cD))=\Ext(\gamma)$, and by \cref{lem:one-boundary}\labelcref{it:one-hull} the sites $y\in\overline\gamma$ with $\dist_*(y,\Ext(\gamma))=1$ are exactly those of $\gamma$; hence $\partialhat B(\cD)=B(\cD)\cap\gamma=\gamma$. Every other component of $\cC$ consists of patterns whose balls are contained in $\Int(\gamma)$, hence its $*$-hull is contained in $\Int(\gamma)$ by \cref{fact:simply}, as $\Int(\gamma)$ is simply connected by \cref{lem:one-boundary}\labelcref{it:one-simply}; so it is nested in $\cD$, and $\cD$ is the unique outermost component of $\cC$. Thus $\cC\in\fC_{\mathfrak c}(\gamma,\eta)$, as required.
\end{enumerate}
\begin{remark}
Part~\labelcref{it:split-map-2} is the point where the connectedness of $\gamma$ and the definition of linked patterns are used; with a boundary that is merely $*$-connected, or with patterns linked only when their balls intersect, the patterns whose balls intersect $\gamma$ could split into several components.
\end{remark}
Having established the bijection, we rewrite \cref{eq:Ec-def} as the following sum over the pairs $(\cC_\gamma,\cC_{\mathrm{in}})$ described above, with $\cC:=\cC_\gamma\cup\cC_{\mathrm{in}}$:
\begin{equation}\label{eq:Ec-pairs}
E_{\mathfrak c}(\gamma,\eta) =\sum_{\cC_\gamma} \sum_{\substack{\cC_{\mathrm{in}}\in\fC(\Int(\gamma))\\ \text{compatible with }\cC_\gamma}} \prod_{\fp\in\cC_\gamma}f(\fp)\prod_{\fp\in\cC_{\mathrm{in}}}f(\fp) \,
\mu_{\Int(\gamma)}^\eta\Big(\sigma\restriction_{B(\cC)\setminus\gamma}=\omega_{\cC}\restriction_{B(\cC)\setminus\gamma}\Big)\,.
\end{equation}
For $\sigma\in \scS^{\Int(\gamma)}$, denote by $\sigma^\eta\in \scS^{\overline\gamma}$ the configuration equal to $\eta$ on $\gamma$ and to $\sigma$ on $\Int(\gamma)$. With this notation, \cref{eq:indicator-product} gives, for every pair in \cref{eq:Ec-pairs},
\[
\mu_{\Int(\gamma)}^\eta\Big(\sigma\restriction_{B(\cC)\setminus\gamma}=\omega_{\cC}\restriction_{B(\cC)\setminus\gamma}\Big) = \E_{\Int(\gamma)}^\eta\bigg[\prod_{\fp\in\cC_\gamma}\one_\fp(\sigma^\eta)\prod_{\fp\in\cC_{\mathrm{in}}}\one_\fp(\sigma)\bigg]\,.
\]
The right-hand side is defined for every $\cC_{\mathrm{in}}\in\fC(\Int(\gamma))$, and vanishes whenever $\cC_{\mathrm{in}}$ is not compatible with $\cC_\gamma$; hence, after this substitution, the sum over $\cC_{\mathrm{in}}$ in \cref{eq:Ec-pairs} may be extended to all of $\fC(\Int(\gamma))$. By \cref{eq:indicator-product,eq:Hg-expansion}, applied with $V=\Int(\gamma)$,
\[
\sum_{\cC_{\mathrm{in}}\in\fC(\Int(\gamma))}\ \prod_{\fp\in\cC_{\mathrm{in}}}f(\fp)\,\one_\fp(\sigma)=e^{\cH_{\Int(\gamma)}^\fg(\sigma)}\,.
\]
Altogether, we arrive at the identity
\begin{equation}\label{eq:Ec-bound}
E_{\mathfrak c}(\gamma,\eta) = \E_{\Int(\gamma)}^\eta\Big[e^{\cH_{\Int(\gamma)}^\fg(\sigma)}\,\Phi_\gamma(\sigma^\eta)\Big]\,,\qquad
\Phi_\gamma(\xi) := \sum_{\cC_\gamma}\ \prod_{\fp\in\cC_\gamma}f(\fp)\,\one_\fp(\xi)\,,
\end{equation}
where in $\Phi_\gamma$ the sum is over the collections $\cC_\gamma\in\fC(\overline\gamma)$ with $\cC_\gamma\subset\cI_\gamma$ and $\gamma\subset B(\cC_\gamma)$ (the constraint $\omega_{\cC_\gamma}\restriction_\gamma=\eta$ is enforced by the indicators $\one_{\fp}(\sigma^\eta)$). By \cref{eq:Ec-bound,eq:Ei-def,eq:w-def},
\begin{equation}\label{eq:w-Phi}
|\sw(\gamma,\eta)| = \bigg|\frac{\E_{\Int(\gamma)}^\eta\big[e^{\cH_{\Int(\gamma)}^\fg(\sigma)}\,\Phi_\gamma(\sigma^\eta)\big]}{\E_{\Int(\gamma)}^\eta\big[e^{\cH_{\Int(\gamma)}^\fg(\sigma)}\big]} \bigg|\leq\|\Phi_\gamma\|_\infty\,.
\end{equation}

Fix $\xi\in \scS^{\overline\gamma}$. A pattern $\fp$ with $\one_\fp(\xi)=1$ is determined by its pair $(x,r)$, as $\omega=\xi\restriction_{B(x,r)}$; thus, if
\[ \cS_\gamma := \big\{(x,r)\,:\; B(x,r)\cap\gamma\neq\emptyset,\, B(x,r)\subset\overline\gamma\big\}\,, \]
then we may rewrite the sum over $\cC_\gamma$ in \cref{eq:Ec-bound} as 
\[
\Phi_\gamma(\xi) = \sum^*_{I} \prod_{(x,r)\in I} f\big(x,r,\xi\restriction_{B(x,r)}\big)\,,
\]
where the $\sum^*$ goes over subsets $I\subset\cS_\gamma$ such that
$\gamma\subset \bigcup_{(x,r)\in I}B(x,r)$.
We claim that every such $I$ in the above sum satisfies
\begin{equation}\label{eq:thinness}
|\gamma|\leq 8\sum_{(x,r)\in I} (r+1)\,.
\end{equation}
Indeed, let $y\in\gamma$. Let $(x,r)\in I$ with $y\in B(x,r)$ (which exists by construction). By \cref{eq:boundary-def} and \cref{lem:one-boundary}\labelcref{it:one-ext}, $y\starsim z$ for some $z\in\Ext(\gamma)$, and $z\notin\overline\gamma\supset B(x,r)$; thus $y$ belongs to the set of sites $u\in B(x,r)$ that are $*$-adjacent to $\Z^2\setminus B(x,r)$, which has at most $8(r+1)$ elements (with room). Summing over $(x,r)\in I$ gives \cref{eq:thinness}. Consequently, by \cref{eq:f-bound}, every such $I$ has
\[
\prod_{(x,r)\in I} \big|f\big(x,r,\xi\restriction_{B(x,r)}\big)\big| \leq e^{-(\Cst-1)\sum_{(x,r)\in I}(r+1)}\leq e^{-\frac{\Cst-1}8|\gamma|}\,,
\]
and in particular, raising this inequality to the power $1-\frac{5}{\Cst-1}>0$ (valid for all $\Cst>6$),
\[
|\Phi_\gamma(\xi)| \leq e^{-\frac{\Cst-6}8|\gamma|}\sum^*_{I}\prod_{(x,r)\in I}\big|f\big(x,r,\xi\restriction_{B(x,r)}\big)\big|^{\frac{5}{\Cst-1}}
\leq e^{-\frac{\Cst-6}8|\gamma|} \sum^*_{I}\prod_{(x,r)\in I} e^{-5(r+1)}\,,
\]
where the last step used \cref{eq:f-bound} again. At this point we can increase $\sum^*_I$ to $\sum_{I\subset\cS_\gamma}$ to get
\begin{align*}
|\Phi_\gamma(\xi)|&\leq e^{-\frac{\Cst-6}8|\gamma|}\prod_{(x,r)\in\cS_\gamma}\Big(1+e^{-5(r+1)}\Big)\\
&\leq \exp\Big(-\tfrac{\Cst-6}{8}|\gamma| + |\gamma|\sum_{r\geq0}(2r+1)^2e^{-5(r+1)}\Big)\\ &\leq e^{-\frac{\Cst-7}8|\gamma|}\,,
\end{align*}
where the second inequality uses that at most $(2r+1)^2|\gamma|$ sites $x$ satisfy $B(x,r)\cap\gamma\neq\emptyset$, and the last one uses $\sum_{r\geq0}(2r+1)^2e^{-5(r+1)}\leq\tfrac18$. For $\Cst\geq 35$ we have $\frac{\Cst-7}8\geq \frac{\Cst}{10}$, and the sought \cref{eq:w-bound} now follows from \cref{eq:w-Phi}. This concludes the proof of \cref{thm:main}, with $C_0=35$.\qed
\end{enumerate}

\subsection{Proof of \texorpdfstring{\cref{cor:potential}}{the effective potential}}\label{subsec:proof-potential}

If $\delta\leq 0$, let $\Psi_A:=g_{x,r}$ when $A=B(x,r)$ for some $(x,r)$ and $\Psi_A:=0$ otherwise; then $\sum_A \Psi_A(\sigma\restriction_A) = \cH_\Lambda^\fg(\sigma)$, and hence the identity in \cref{eq:effective-potential} is nothing but \cref{eq:HD-def}, while $\sum_{A\ni y}\|\Psi_A\|_\infty e^{\delta|A|}\leq\sum_{r\geq0}(2r+1)^2e^{-\Cst (r+1)}\leq 1 \leq e^{-\delta}$. Assume henceforth that $\delta>0$. 

We will obtain the corollary from \cref{thm:main} by a cluster expansion performed for each fixed~$\sigma$. (The base measure $\mu_\Lambda$ plays no role in this step beyond defining the weights $\sw$.)
To do so, recall the following standard notation (see, e.g., \cite[\S5]{FriedliVelenik18}). The Ursell functions in \cref{eq:cluster_expansion} are given by
\begin{equation}\label{eq:ursell-def}
\phi_\zeta(\pi_1,\ldots,\pi_n):=\sum_{G}\ \prod_{ij\in E(G)}\big(\zeta(\pi_i,\pi_j)-1\big)\,,
\end{equation}
the sum going over connected graphs $G$ on the vertices $\{1,\ldots,n\}$, with $E(G)$ denoting its edge~set. For a finite set of polymers with a symmetric compatibility relation $\perp$ such that $\pi\not\perp\pi$, and complex activities $\sw$, we take $\zeta(\pi_1,\pi_2)=\one_{\{\pi_1\perp\pi_2\}}$, and denote the left-hand side of \cref{eq:cluster_expansion} by $\Xi:=\sum_{\cF}\prod_{\pi\in\cF}\sw(\pi)$ (the sum is over families $\cF$ of pairwise compatible polymers). Further denote the term of order $n$ in the sum in \cref{eq:cluster_expansion}, and the corresponding sum of absolute values, by
\begin{equation}\label{eq:cn-def}
\fb_n:=\frac1{n!}\sum_{(\pi_1,\ldots,\pi_n)}\phi_\zeta(\pi_1,\ldots,\pi_n)\prod_{i=1}^n\sw(\pi_i)\qquad\mbox{and}\qquad
\overline{\fb}_n:=\frac1{n!}\sum_{(\pi_1,\ldots,\pi_n)}\big|\phi_\zeta(\pi_1,\ldots,\pi_n)\big|\prod_{i=1}^n|\sw(\pi_i)|\,,
\end{equation}
respectively, the sums going over the ordered $n$-tuples of polymers. Whenever \cref{eq:textbook_bound} holds, one has $\sum_n\overline{\fb}_n<\infty$, $\log\Xi=\sum_n\fb_n$, and for every polymer $\pi_0$,
\begin{equation}\label{eq:cluster-bound-general}
\sum_{n\geq1}\frac1{n!}\sum_{\substack{(\pi_1,\ldots,\pi_n)\\ \exists i\,:\,\pi_i\not\perp\pi_0}}\big|\phi_\zeta(\pi_1,\ldots,\pi_n)\big|\prod_{i=1}^n |\sw(\pi_i)|\leq a(\pi_0)\,,
\end{equation}
by \cite[Thm.~5.4 and Eqs.~(5.10)--(5.11)]{FriedliVelenik18} (multiply the latter by $|\sw(\pi_1)|$ and sum over $\pi_1\not\perp\pi_0$, the factor $k$ there accounting for the choice of the index $i$, as in \cite[Eq.~(5.29)]{FriedliVelenik18}).

Fix $\sigma\in \scS^\Lambda$. In \cref{eq:polymer-rep}, for a given support $\gamma$, only the polymer $(\gamma,\sigma\restriction_\gamma)$ contributes, whence \begin{equation}\label{eq:Z-Lambda-via-Xi}
Z_\Lambda^\fg = \E[\Xi(\sigma)] \qquad\mbox{for}\qquad \Xi(\sigma) = \sum_{\cF\subset\Gamma}\ \prod_{\gamma\in\cF}\sw(\gamma,\sigma\restriction_\gamma)\,,
\end{equation}
the sum running over families of polymers with pairwise compatible supports. Put 
\[\widetilde{\sw}(\gamma):=e^{2\delta|\gamma|}\sup_{\eta}|\sw(\gamma,\eta)|\,,\] 
so that $|\sw(\gamma,\sigma\restriction_\gamma)|\leq e^{-2\delta|\gamma|}\,\widetilde{\sw}(\gamma)$ and, by \cref{eq:w-bound}, 
\[ \widetilde{\sw}(\gamma)\leq e^{-(\Cst/10-2\delta)|\gamma|}=e^{-4|\gamma|}\,.\] 
We will now verify the criterion~\labelcref{eq:textbook_bound} for the activities $\widetilde{\sw}$ with $\zeta(\gamma_1,\gamma_2)=\one_{\{\gamma_1\perp\gamma_2\}}$ and $a(\gamma)=|\gamma|$. Fix $\gamma_1\in\Gamma$. If $\gamma_2\not\perp\gamma_1$ then $\dist(\gamma_1,\gamma_2)\leq1$ by \cref{eq:compatible-def}, so $\gamma_2$ contains one of the at most $5|\gamma_1|$ sites $z$ with $\dist(z,\gamma_1)\leq 1$. Every $\gamma_2\in\Gamma$ is connected by \cref{lem:one-boundary}\labelcref{it:one-conn}, and the number of connected sets of size $n$ rooted at a given site in a graph of maximum degree $d$ is at most the number of subtrees with $n$ vertices, containing the root, of the infinite $d$-ary tree, which is the Fuss--Catalan number $\frac1{(d-1)n+1}\binom{dn}n\leq(ed)^{n-1}$ (see, e.g., \cite[p.~397, Ex.~11]{Knuth_TAOCP1}); here we have $d=4$. 
Hence,
\[
\sum_{\gamma_2\,:\,\gamma_2\not\perp\gamma_1} \widetilde{\sw}(\gamma_2)\,e^{|\gamma_2|}\leq 5|\gamma_1|\sum_{n\geq1}(4e)^{n-1} e^{-3n} = \frac{5}{e^3-4e}|\gamma_1|<|\gamma_1|\,,
\]
and so \cref{eq:textbook_bound} holds for $\widetilde\sw(\gamma)$, and thus also for $\sw(\gamma,\sigma\restriction_\gamma)$ (since $|\sw(\gamma,\sigma\restriction_\gamma)|\leq\widetilde{\sw}(\gamma)$ as $\delta\geq 0$). Hence, $\log\Xi(\sigma)=\sum_n\fb_n$ with $\fb_n$ as in \cref{eq:cn-def} for the activities $\sw(\gamma,\sigma\restriction_\gamma)$, and grouping the $n$-tuples there according to the set $A:=\gamma_1\cup\cdots\cup\gamma_n$ shows that $\Xi(\sigma) = \exp\big(\sum_A \Psi_A(\sigma\restriction_A)\big)$ for
\begin{equation}\label{eq:Psi-def}
 \Psi_A(\sigma) := \sum_{n\geq1}\frac1{n!}\sum_{\substack{(\gamma_1,\ldots,\gamma_n)\\ \gamma_1\cup\cdots\cup\gamma_n=A}}\phi_\zeta(\gamma_1,\ldots,\gamma_n)\prod_{i=1}^n\sw(\gamma_i,\sigma\restriction_{\gamma_i})\,.
\end{equation}
Note that $\Psi_A$ depends on $\sigma$ only through $\sigma\restriction_A$, and only on $\fg$ and the measures $\mu_{\Int(\gamma_i)}^\eta$ through the weights. Whenever $\phi_\zeta(\gamma_1,\ldots,\gamma_n)\neq0$ the incompatibility graph is connected; since incompatible supports are at distance at most $1$ and each support is connected, $A$ is then connected, so $\Psi_A=0$ unless $A$ is connected. Combining this with \cref{eq:Z-Lambda-via-Xi}, we arrive at the identity in \cref{eq:effective-potential}. 

For the inequality in \cref{eq:effective-potential}, fix $y\in\Lambda$; using that $2\delta\sum |\gamma_i| \geq \delta|A| + \delta |\gamma_1| \geq \delta|A| + \delta$, we get
\[
\sum_{A\ni y}\|\Psi_A\|_\infty e^{\delta|A|}\leq \sum_{n\geq1}\frac1{n!}\sum_{\substack{(\gamma_1,\ldots,\gamma_n)\\ y\in\gamma_1\cup\cdots\cup\gamma_n}}\big|\phi_\zeta(\gamma_1,\ldots,\gamma_n)\big|e^{-\delta}\prod_{i=1}^n \widetilde{\sw}(\gamma_i)\leq e^{-\delta}\, a(\{y\}) = e^{-\delta}\,,
\]
where the second inequality used \cref{eq:cluster-bound-general} for the activities $\widetilde{\sw}(\gamma)$ (which satisfy \cref{eq:textbook_bound}) and the polymer $\gamma_0:=\{y\}=\partialhat\{y\}\in\Gamma$ (any $\gamma\ni y$ is incompatible with $\gamma_0$).
This completes the proof.
\qed

\section{Cluster expansions: remarks and applications}\label{sec:corollaries}

Throughout this section, we say that the potentials $\scV=(\scV_{x,y},\scV_x)$ are defined in $\Z^2$ if they are defined for all sites $x$ and edges $xy$ in $\Z^2$ independently of the domain $\Lambda\subset\Z^2$, and that they are translation invariant if they do not depend on $x$; similarly, for the tilt $\fg=(g_{x,r})$.

\begin{proposition}[Product measure]\label{prop:product}
In the setting of \cref{thm:main}, with $\mu_\Lambda$ a product measure ($\scV_{x,y}\equiv 0$ in \labelcref{eq:gibbs-def}), there are absolute constants $C,C',c>0$ so that for $\Cst\geq C$ the following holds.
\begin{enumerate}[(a)]
\item\label{it:prod-KP} The series in \cref{eq:cluster_expansion} for the polymers and the compatibility relation of \cref{thm:main}, with the activities $w(\gamma,\eta):=\sw(\gamma,\eta)\,\mu_\Lambda(\sigma\restriction_\gamma=\eta)$, converges absolutely to $\log Z_\Lambda^{\fg}$.
\item\label{it:prod-local} For every finite $X\subset\Lambda$ and $\xi\in\scS^{X}$ with $\mu_\Lambda(\sigma\restriction_{X}=\xi)>0$, and every finite $Y\subset\Lambda$ disjoint from $X$ and $\xi'\in\scS^{Y}$ with $\mu_\Lambda(\sigma\restriction_{Y}=\xi')>0$,
\begin{equation}\label{eq:prod-local-3}
\Big|\log\frac{\mu_\Lambda^{\fg}(\sigma\restriction_{X}=\xi,\ \sigma\restriction_{Y}=\xi')}{\mu_\Lambda^{\fg}(\sigma\restriction_{X}=\xi)\,\mu_\Lambda^{\fg}(\sigma\restriction_{Y}=\xi')}\Big|\leq C'|X|\,e^{-c\,\dist(X,Y)}\,.
\end{equation}
Further suppose $\scV,\fg$ are defined in $\Z^2$. Then, for every finite simply connected $\Lambda'\supset\Lambda$,
\begin{equation}\label{eq:prod-local-2}
\Big|\log\frac{\mu_{\Lambda'}^{\fg}(\sigma\restriction_{X}=\xi)}{\mu_\Lambda^{\fg}(\sigma\restriction_{X}=\xi)}\Big|\leq C'|X|\,e^{-c\,\dist(X,\Z^2\setminus\Lambda)}\,.
\end{equation}
Consequently $\mu_\Lambda^{\fg}$ converges weakly, as $\Lambda\uparrow\Z^2$, to a probability measure $\mu_{\infty}^{\fg}$ on $\scS^{\Z^2}$, which satisfies \cref{eq:prod-local-3}, and which is translation invariant if $\scV$ and $\fg$ are.
\end{enumerate}
\end{proposition}

\begin{remark}\label{rem:boundary-tilt}
The hypothesis that $\fg$ is defined in $\Z^2$ in \cref{prop:product}\labelcref{it:prod-local}
can be relaxed: each finite $\Lambda$ may modify the terms $g_{x,r}$ with $B(x,r)\cap\partialvtxin\Lambda\neq\emptyset$, subject to \cref{eq:g-decay}, without affecting \cref{eq:prod-local-2} or the limit $\mu_\infty^\fg$. For instance, a boundary condition may be encoded in the terms $(g_{x,0})_{x\in\partialvtxin\Lambda}$.
\end{remark}

\begin{remark}\label{rem:highT}
\Cref{prop:product} covers any nearest-neighbor model at sufficiently high temperature, as a weak interaction of range $\rho$ can be absorbed into $g_{x,\rho}$: for the Ising model of \cref{eq:ising-def}, take $\mu_\Lambda$ to be the uniform measure, replace $g_{x,1}$ by $g_{x,1}-\frac\beta2\sum_{y\in\Lambda\,:\,y\sim x}\one_{\{\sigma_x\neq\sigma_y\}}$, and encode the boundary condition in $(g_{x,0})_{x\in\partialvtxin\Lambda}$ (via \cref{rem:boundary-tilt}). The new tilt satisfies \cref{eq:g-decay} once $\beta$ is small enough; the cost is that the admissible range of $\beta$ depends on $\Cst$, unlike the absolute bound on $\beta$ obtained by expanding the base model and the effective potential of \cref{cor:potential} separately.
\end{remark}

\begin{proposition}[Low temperature Ising]\label{prop:lowT}
In the setting of \cref{thm:main}, with $\Z^2\setminus\Lambda$ connected, $\mu_\Lambda$ the Ising model~\labelcref{eq:ising-def} with all-plus boundary conditions ($\scV_{x,y}=\beta\one_{\{s\neq s'\}}$, $\scV_x\equiv0$, and $\tau=\underline{1}$ (the all-plus configuration) in \labelcref{eq:gibbs-def}), and with $\fg$ spin-flip symmetric ($g_{x,r}(-\omega)=g_{x,r}(\omega)$), there are absolute constants $C,\beta_0,C',c>0$ so that for $\Cst\geq C$ and $\beta\geq\beta_0$ the following holds.
\begin{enumerate}[(a)]
\item\label{it:ising-KP} There are activities $W(S)$, indexed by $*$-connected sets $S$ of dual sites $z\in \Z^2+(\frac12,\frac12)$ with $\dist_*(z,\Lambda)\leq\frac12$ and compatible when disjoint, for which the series $\sum_n\fb_n$ in \cref{eq:cn-def} converges absolutely to $\log(Z_\Lambda^\tau Z_\Lambda^{\fg})-\sum_A\Psi_A(\underline{1})$, where $\Psi_A$ is the potential of \cref{cor:potential}.
\item\label{it:ising-local} For all finite disjoint subsets $X,Y\subset\Lambda$, 
\begin{equation}
\sup_{\xi,\xi'}\Big|\mu_\Lambda^{\fg}(\sigma\restriction_{X}=\xi,\ \sigma\restriction_{Y}=\xi')-\mu_\Lambda^{\fg}(\sigma\restriction_{X}=\xi)\,\mu_\Lambda^{\fg}(\sigma\restriction_{Y}=\xi')\Big|\leq e^{C'(|X|+|Y|)}\,e^{-c\,\dist(X,Y)}\,.\label{eq:ising-local-3}
\end{equation}
Further suppose $\fg$ is defined in $\Z^2$. Then, for every finite $\Lambda'\supset\Lambda$ with $\Z^2\setminus\Lambda'$ connected,
\begin{equation}
\sup_{\xi}\Big|\mu_{\Lambda'}^{\fg}(\sigma\restriction_{X}=\xi)-\mu_\Lambda^{\fg}(\sigma\restriction_{X}=\xi)\Big|\leq e^{C'|X|}\,e^{-c\,\dist(X,\Z^2\setminus\Lambda)}\,.\label{eq:ising-local-2}
\end{equation}
Consequently $\mu_\Lambda^{\fg}$ converges weakly, as $\Lambda\uparrow\Z^2$, to a probability measure $\mu_{\infty}^{\fg}$ on $\{\pm1\}^{\Z^2}$, which satisfies \cref{eq:ising-local-3} and is translation invariant if $\fg$ is.
\end{enumerate}
\end{proposition}
The spin-flip symmetry of $\fg$ in \cref{prop:lowT} guarantees that it does not break the symmetry between the two low-temperature phases, so that the Peierls contours of the base model remain the relevant polymers; without it, the tilt acts as a (weak, but volume-extensive) field selecting a phase, and the appropriate framework is that of Pirogov--Sinai theory, which we do not enter here.

We prove these in \cref{subsec:proof-product,subsec:proof-lowT}, after discussing when the classical expansion might fail and the analyticity of $\log Z_\Lambda^\fg$ (\cref{subsec:direct-nonanalytic}), and formulating a polymer comparison tool (\cref{subsec:comparison}).

\subsection{Examples of the classical expansion failing and non-analyticity of the free energy}\label{subsec:direct-nonanalytic}
\subsubsection{An example where the classical expansion fails}\label{subsec:direct}

We show that for the classical expansion in \cref{ex:product-abs-mag}, whose polymers are obtained by grouping the balls in \labelcref{eq:classical-expansion} into connected unions, the series in \labelcref{eq:cluster_expansion} may diverge already in finite volume, even though the model satisfies the assumptions of \cref{thm:main} and \cref{prop:product}. Call a collection of balls connected if the graph on the collection in which two balls are adjacent when they intersect is connected; with an i.i.d.\ base measure, and with each $S_i$ below the union of a connected collection of balls $B(x_j,r_j)$, \cref{eq:classical-expansion} reads
\begin{equation}\label{eq:w(S)-naive}
 Z_\Lambda^{\fg} = \sum_{\substack{\{S_i \}\\ \text{pairwise disjoint}}} \prod_i \sw(S_i)\qquad\mbox{for}\qquad \sw(S):= \sum_{\substack{\{ (x_{j}, r_{j}) \} \text{ s.t.}\\ \bigcup_j B(x_{j}, r_{j}) = S} }  \E \bigg[\prod_j \Big(e^{ g_{x_{j}, r_{j}} ( \sigma \restriction_{B(x_{j}, r_{j})}) } -1\Big)\bigg]\,,
\end{equation}
and \cref{eq:cluster_expansion} with the interaction $\zeta(S_1,S_2)=\one_{\{S_1\cap S_2=\emptyset\}}$ then reads
\[
 Z_\Lambda^{\fg}  = \sum_{ \{S_i \}} \prod_i \sw( S_i) \prod_{i, j} \one_{\{ S_i \cap  S_j = \emptyset\}}  = \exp\bigg[ \sum_{n\geq1}\frac1{n!} \sum_{(S_1, \ldots, S_n)} \phi_\zeta( S_1, \ldots, S_n)\prod_i \sw(S_i) \bigg]\,.
\]

\begin{proposition}[Divergence of the classical expansion]\label{prop:naive}
For $\mu_\Lambda^\fg$ as in \cref{ex:product-abs-mag} and every constant $C>0$, there exist $N_0=N_0(C)$ and $\kappa=\kappa(C)>1$ such that the series in \cref{eq:cluster_expansion} for the classical expansion in an $N\times N$  box diverges for any $N\geq N_0$. More precisely, the quantities $\fb_n$ and $\overline{\fb}_n$ as defined in \cref{eq:cn-def} (so that the right-hand side of \labelcref{eq:cluster_expansion} is $\exp(\sum_n\fb_n)$) satisfy
\[
(-1)^{n-1}\fb_n=\overline{\fb}_n\geq \kappa^n / n \qquad\mbox{for all $n\geq1$}\,.
\]
\end{proposition}
\noindent (NB.\ \Cref{thm:main,prop:product} apply to this model, uniformly in $N$, if $C$ is large enough.)
\begin{proof}
We will need the next result, which follows from classical properties of the Tutte polynomial of multigraphs; see \cref{app:markov} for references and a short direct proof, provided for completeness.
\begin{claim}\label{clm:ursell-sign}
Let $\zeta(\pi_1,\pi_2) = \one_{\{\pi_1 \perp \pi_2\}}$. Then one has $(-1)^{n-1}\phi_\zeta(\pi_1,\ldots,\pi_n)\geq0$ for every $n\geq1$ and $\pi_1,\ldots,\pi_n$; moreover, if $\pi_1=\cdots=\pi_n$ then $\phi_\zeta(\pi_1,\ldots,\pi_n)=(-1)^{n-1}(n-1)!$.
\end{claim}
Since $|B(x,r)|=(2r+1)^2$ is odd, $|\sum_{v\in B(x,r)}\sigma_v|\geq1$, and hence
\[
e^{g_{x,r}}-1\geq g_{x,r}\geq e^{-C(r+1)}\qquad\mbox{pointwise, for all $x$ and $r$}\,.
\]
In particular, all the summands in the definition of $\sw(S)$ in \cref{eq:w(S)-naive} are nonnegative, implying that $\sw(S)\geq 0$ for every polymer $S$, and therefore $(-1)^{n-1}\fb_n=\overline{\fb}_n$ for every $n$ by \cref{clm:ursell-sign}. By the second part of that claim, it suffices to find a single polymer $S_0$ with $\sw(S_0) \geq \kappa$ for some $\kappa>1$: 
taking $S_i=S_0$ for $i=1,\ldots,n$, this gives $|\phi_\zeta(S_0,\ldots,S_0)|=(n-1)!$, and hence 
\[ \overline{\fb}_n\geq\frac1{n!}(n-1)!\,\sw(S_0)^n\geq\frac{\kappa^n}n\,.\] 
To that end, let $S_0:=B(x_0,r_0)\subset\Lambda$ for some large enough $r_0=r_0(C)$. Since every term $e^{g_{x,r}}-1$ in the definition of $\sw(S)$ in \labelcref{eq:w(S)-naive} is nonnegative, we may bound $\sw(S_0)$ from below via the full $S_0$ (a $g_{x_0,r_0}$ term) and every set $S\subset S_0$ of singletons $x$ (with $e^{g_{x,0}}-1=e^{e^{-C}}-1$ terms):
\begin{align*}
\sw(S_0) &\geq \sum_{S\subset S_0} \E\Big[ g_{x_0,r_0}(\sigma\restriction_{S_0})  \prod_{x\in S} \Big(e^{e^{-C}}-1\Big)\Big] 
= \E \Big[ g_{x_0,r_0}(\sigma\restriction_{S_0}) \prod_{x\in S_0} e^{e^{-C}}\Big] \\
&\geq \exp\Big( (2r_0+1)^2 e^{-C}-C(r_0+1)\Big) =: \kappa\,,
\end{align*}
which exceeds $1$ provided that $r_0$ is large enough in terms of $C$. Fixing such an $r_0$, let $N_0:=2r_0+1$, so that for $N\geq N_0$ the box $\Lambda$ indeed contains a ball $S_0$ as above, which concludes the proof.
\end{proof}

\begin{remark}[Shifting constants]\label{rem:direct-shift}
The proof of \cref{prop:naive} exhibits a single polymer of weight larger than $1$, and it uses the pointwise lower bounds $e^{g_{x,r}}-1\geq e^{-C(r+1)}$. Replacing each $g_{x,r}$ by $g_{x,r}-c_{x,r}$ for constants $c_{x,r}$ multiplies $Z_\Lambda^\fg$ by a constant and changes the weights of the classical expansion, and the proof above does not apply to the shifted expansion. We do not expect such a shift to restore absolute convergence: let us argue that a shift cannot remove the small polymers.
For $x\sim y$ with $B(x,1)\cup B(y,1)\subset \Lambda$---so $\E[e^{g_{x,1}}]=\E[e^{g_{y,1}}]=:\upmu$ as the base measure is i.i.d.---let
\[
S_{x,y}:=B(x,1)\cup B(y,1)\qquad\mbox{and}\qquad \uprho  :=\Cov\big(e^{g_{x,1}} , e^{g_{y,1}}\big)/\upmu^2\,.
\]
We claim $\uprho>0$. 
If $\fm = \sum_{z\in B(x,1)\cap B(y,1)}\sigma_z$ then $f_x(m):=\E [e^{g_{x,1}} \mid \fm =m]$ is even and convex (namely, $\E[\exp(e^{-2C}|m+Y|)]$ where $Y:=\sum\{\sigma_z:z\in B(x,1)\setminus B(y,1)\}\stackrel{\mathrm{d}}{=}-Y$ is independent of $\fm$), thus nondecreasing in $|m|$, and similarly for $f_y(m)$.
Each is nonconstant, so $\Cov(f_x(\fm),f_y(\fm))>0$ by Chebyshev's sum inequality. The conditional independence of $e^{g_{x,1}}, e^{g_{y,1}}$ given $\fm$ implies that $\E[e^{g_{x,1}+g_{y,1}}]=\E[f_x(\fm)f_y(\fm)]$; thus, $\Cov(e^{g_{x,1}},e^{g_{y,1}})=\Cov(f_x(\fm),f_y(\fm))>0$ and so $\uprho>0$.

Now let $\hat\sw$ be the weights in \cref{eq:w(S)-naive} for $\hat g_{x,r} =g_{x,r}-\Psi_r$ for all $x$, for some $\Psi_r$, $r=0,1$. Then
\[
\hat\sw\big(S_{x,y}\big)=e^{12(e^{-C}-\Psi_0)}\Big(
  e^{-2\Psi_1} (1+\uprho) \upmu^2 - 2e^{-\Psi_1} \upmu + 1\Big) \geq e^{12(e^{-C}-\Psi_0)}\,\frac{\uprho}{1+\uprho}\,,
\]
the prefactor arising from the $12$ singletons $\{v\}\subset S_{x,y}$, each of activity $e^{e^{-C}-\Psi_0}-1$, an arbitrary set of which may be added to the pair $\{B(x,1),B(y,1)\}$ (a collection with this union that contains, e.g., $B(x,1)$ but not $B(y,1)$, cannot be connected in the intersecting sense); and $1+\uprho>0$. 

Finally, if $|\Psi_0-e^{-C}|\geq1$ then every singleton $\{x\}$ has $|\hat\sw(\{x\})| = |\exp(e^{-C}-\Psi_0)-1|\geq 1-e^{-1}$; otherwise, the pairs $x\sim y$ have $\hat\sw(S_{x,y}) \geq e^{-12}\uprho/(1+\uprho)$. Either way,~$\Lambda$ has $\Omega(N^2)$ pairwise disjoint polymers with weights bounded away from $0$, uniformly in $\Psi_0,\Psi_1$ and~$N$.
\end{remark}

\subsubsection{Non-analyticity of the free energy}\label{subsec:analyticity}

As mentioned in \cref{rem:analyticity}, we now give a heuristic, based on cumulants, suggesting that the free energy $\lim_\Lambda|\Lambda|^{-1}\log Z_\Lambda^{\fg}$ should not, in general, be analytic in a global prefactor of the interactions under \cref{eq:g-decay}, unlike in classical cluster expansion.

\begin{remark}\label{rem:analyticity-obstruction}
The following candidate should have a non-analytic free energy as a function of the interaction. Let $\mu_\Lambda$ be the uniform measure on $\{\pm1\}^\Lambda$ and consider $\fg=\beta\fg^0$, so that $\partial_\beta^k\log Z_\Lambda^{\fg}|_{\beta=0}$ is the $k$-th cumulant $\kappa_k(\cH^0)$ of $\cH^0=\sum_{x,r}g^0_{x,r}$ under $\mu_\Lambda$; analyticity in a $\Lambda$-independent complex neighborhood of $\beta=0$ would force $|\kappa_k(\cH^0)|\leq|\Lambda|\,k!\,R^{-k}$ for some $R>0$. Take 
\begin{align*}
g^0_{x,0}(\omega) &:=e^{-\Cst}\omega_x\,,\\
g^0_{x,r}(\omega) &:=e^{-\Cst (r+1)}\sign P_{x,r}(\omega)\qquad(r\geq 1)\,,\ \end{align*} 
for
\[
P_{x,r}(\omega):=\sum_{\{y_1,\ldots,y_{\lceil \Cst r\rceil-1}\}\subset B(x,r-1)}\prod_{i=1}^{\lceil \Cst r\rceil-1}\omega_{y_i}\,,
\]
where $P_{x,r}$ is the elementary symmetric polynomial of degree $\lceil \Cst r\rceil-1$ in the spins of the ball. Fix $k\geq1$ and let $r$ be such that $\lceil \Cst r\rceil=k$, so that $r\approx k/\Cst$. By multilinearity, $\kappa_k(\cH^0)$ is the sum over ordered $k$-tuples of balls of the joint cumulants of the corresponding $g^0$'s. Since the spins are independent and centered, the joint cumulant of $g^0_{x,r}$ with $g^0_{y_1,0},\ldots,g^0_{y_{k-1},0}$, for distinct $y_i\in B(x,r-1)$, equals $e^{-\Cst(k-1)}\E_\Lambda[g^0_{x,r}\prod_i\omega_{y_i}]$, so the contribution to $\kappa_k(\cH^0)$ of the tuples consisting of $g^0_{x,r}$ and $k-1$ such balls of radius zero is, summing over the ordered $(k-1)$-tuples $(y_1,\ldots,y_{k-1})$ and using that $\sum_{\{y_i\}}\prod_i\omega_{y_i}=P_{x,r}$, exactly \[(k-1)!\, e^{-\Cst(k-1)}e^{-\Cst (r+1)}\E_\Lambda\big[\sign(P_{x,r})P_{x,r}\big]=(k-1)!\, e^{-\Cst(k-1)}e^{-\Cst (r+1)}\E_\Lambda|P_{x,r}|\,.\] By the Bonami--Beckner inequality, 
\[ \E_\Lambda|P_{x,r}|\geq3^{-(k-1)}\|P_{x,r}\|_2=3^{-(k-1)}\binom{|B(x,r-1)|}{k-1}^{1/2}\,,\]
 which is at least $(c\sqrt k/\Cst)^{k}$ for an absolute $c>0$. Thus, summing over the $|\Lambda|$ positions $x$ of the ball of radius $r$, this contribution to $\kappa_k(\cH^0)$ is of order $|\Lambda|\,(k!)^{3/2}R_{\Cst}^{-k}$, with $R_{\Cst}$ depending on $\Cst$ only, which is not of the form $|\Lambda|\,k!\,R^{-k}$. 
 
 We have not verified that other tuples do not cancel this contribution, but we see no mechanism for such a cancellation; if there is none, then the free energy of this model is not analytic at $\beta=0$, although \cref{prop:product} controls $\log Z_\Lambda^{\fg}$ for all real $\beta$ with $|\beta|\leq1$. 

This is the mechanism of Griffiths singularities\footnote{The example given in~\cite{Griffiths69} is the Ising model on the graph obtained from supercritical bond percolation on $\Z^d$ with parameter $p$, at $\beta$ above the critical point $\beta_c(\Z^d)$ but below the critical point of the diluted model $\beta_c((\Z^d)_p)$; its (quenched) free energy $\lim_{|\Lambda|\to\infty}\E_p [\frac1{|\Lambda|} \log Z_{\Lambda}]$ is not analytic in the external field $h$ at $h=0$.}: a rare structure (here a ball of radius $r$, of weight $e^{-\Cst (r+1)}$) contributes to the $k$-th cumulant $\kappa_k$ some $\Omega(r^{f(k)})$ term where $f(k)$ diverges with~$k$ (here of order $r^{2(k-1)}$, the number of ways to place the $k-1$ balls of radius zero inside it), so that summing over $r$ yields super-factorial growth in $k$. Under the decay $e^{-\Cst r^2}$ of \cite{GMT17,GMT20} the same computation gives $\sup_r e^{-\Cst r^2}r^{k-1}\asymp(k!)^{1/2}$ in place of $\sup_r e^{-\Cst r}r^{k-1}\asymp k!$, hence contributions of order $k!R^{-k}$, consistent with the analyticity obtained there from the direct expansion. Analyticity in a single $g_{x,r}$, or in finitely many, is trivial: $Z_\Lambda^{\fg}$ is then a finite positive combination of exponentials of the variables, nonvanishing as long as their imaginary parts sum to less than $\pi/2$.
\end{remark}

\subsection{Comparison of polymer models}\label{subsec:comparison}
The following lemma allows one to compare a polymer model of a base measure $\mu$ with that of $\mu(\cdot\mid\sigma\restriction_X=\xi)$ for a finite $X$, as well as in a larger domain.

\begin{lemma}\label{lem:comparison}
Let $\fP$ be a finite set of polymers with a symmetric relation $\perp$ such that $\pi\not\perp\pi$; equip each $\pi\in\fP$ with a size $|\pi|\geq1$ and a finite \emph{shadow} $\shadow(\pi)\subset\Z^2$ such that $\diam_*(\shadow(\pi))\leq10|\pi|$ and $\shadow(\pi)\cap\shadow(\pi')\neq\emptyset$ whenever $\pi\not\perp\pi'$. Let $w^{(i)}$ ($i=1,\ldots,4$) be complex activities on $\fP$, and $\Sigma_i:=\sum_{n\geq1}\fb_n^{(i)}$ with $\fb_n^{(i)}$ as in \labelcref{eq:cn-def}.
Suppose that there exist $K,s>0$ so that, for each $i$,
\begin{equation}\label{eq:cmp-KP}
\sum_{\pi'\in\fP\,:\,\pi'\not\perp\pi}|w^{(i)}(\pi')|\,e^{(1+s)|\pi'|}\leq|\pi|\qquad\mbox{for all $\pi\in\fP$}
\end{equation}
and
\begin{equation}\label{eq:shadow-hyp}
\sum_{\pi\in\fP\,:\,\shadow(\pi)\ni y}|w^{(i)}(\pi)|\,e^{(1+s)|\pi|}\leq K\qquad\mbox{for all $y\in\Z^2$}\,.
\end{equation}
The activities $w,w'$ are said to \emph{agree outside} a given set $U\subset\Z^2$, denoted $w\agreeoutside{U^c}w'$, if they coincide on all $\pi$ with $\shadow(\pi)\cap U=\emptyset$. Then, whenever  
$ w^{(1)} \agreeoutside{U^c} w^{(2)}$, we have
\begin{equation}\label{eq:log-X1-X2}
|\Sigma_1-\Sigma_2|\leq 2K|U| \,,
\end{equation}
and if, in addition, $ w^{(3)} \agreeoutside{U^c} w^{(4)}$, $ w^{(1)} \agreeoutside{V^c} w^{(3)}$, $ w^{(2)} \agreeoutside{V^c} w^{(4)}$,  then 
\begin{equation}\label{eq:log-X1-X2-log-X3-X4} \big|(\Sigma_1-\Sigma_2)-(\Sigma_3-\Sigma_4)\big|\leq 4K|U|\,e^{-s\,\dist_*(U,V)/10}\,.
\end{equation}
\end{lemma}

\begin{proof}
Each $w^{(i)}$ satisfies \cref{eq:textbook_bound} with $a(\pi)=|\pi|$ by \cref{eq:cmp-KP}, so the series $\Sigma_i$ converge absolutely. For a polymer $\pi$ and activities $w$ satisfying \cref{eq:textbook_bound} with $a(\pi)=|\pi|$, i.e., \cref{eq:cmp-KP} with $s=0$, we have
\begin{equation}\label{eq:cluster-through}
\sum_{n\geq1}\frac1{n!}\sum_{\substack{(\pi_1,\ldots,\pi_n)\\ \pi\in\{\pi_1,\ldots,\pi_n\}}}\big|\phi_\zeta(\pi_1,\ldots,\pi_n)\big|\prod_{i=1}^n|w(\pi_i)|\leq|w(\pi)|\,e^{|\pi|}\,,
\end{equation}
by \cite[Eq.~(5.11)]{FriedliVelenik18}, as in the derivation of \cref{eq:cluster-bound-general}. For $\fX=(\pi_1,\ldots,\pi_n)$ we write $w(\fX):=\prod_i w(\pi_i)$ and $\shadow(\fX):=\bigcup_i\shadow(\pi_i)$; thus, if $w\agreeoutside{U^c} w'$ then $w(\fX)=w'(\fX)$ when $\shadow(\fX)\cap U=\emptyset$.
With this notation, by \cref{eq:cn-def}, \[\Sigma_1-\Sigma_2=\sum_n\frac1{n!}\sum_\fX\phi_\zeta(\fX)[w^{(1)}(\fX)-w^{(2)}(\fX)]\,,\] and the terms of this series vanish unless $\shadow(\fX)$ intersects $U$. By \cref{eq:cluster-through}, for activities $w$ satisfying \cref{eq:textbook_bound} with $a(\pi)=|\pi|$ and $\sum_{\pi\,:\,\shadow(\pi)\ni y}|w(\pi)|e^{|\pi|}\leq K$ for all $y\in\Z^2$ (as each $w^{(i)}$ does, by \cref{eq:cmp-KP,eq:shadow-hyp} for $s=0$),
\begin{equation}\label{eq:clusters-meeting-U}
\sum_{n\geq1}\frac1{n!}\sum_{\fX\,:\,\shadow(\fX)\cap U\neq\emptyset}\big|\phi_\zeta(\fX)\big|\,|w(\fX)|\leq\sum_{\pi\,:\,\shadow(\pi)\cap U\neq\emptyset}|w(\pi)|\,e^{|\pi|}\leq K|U|\,,
\end{equation}
which, once applied to $w=w^{(1)}$ and $w=w^{(2)}$, gives \cref{eq:log-X1-X2}.

Next, the combination $w^{(1)}(\fX)-w^{(2)}(\fX)-w^{(3)}(\fX)+w^{(4)}(\fX)$ vanishes unless $\shadow(\fX)$ intersects both $U$ and $V$: if $\shadow(\fX)\cap V=\emptyset$ then $w^{(1)}(\fX)=w^{(3)}(\fX)$ and $w^{(2)}(\fX)=w^{(4)}(\fX)$, and if $\shadow(\fX)\cap U=\emptyset$ then $w^{(1)}(\fX)=w^{(2)}(\fX)$ and $w^{(3)}(\fX)=w^{(4)}(\fX)$. Let $\fX=(\pi_1,\ldots,\pi_n)$ be such that $\shadow(\fX)$ intersects both, with $\phi_\zeta(\fX)\neq0$, so that the graph $G(\pi_1,\ldots,\pi_n)$ is connected, and let $\shadow(\pi_a)\cap U\neq\emptyset$ and $\shadow(\pi_b)\cap V\neq\emptyset$. Along a path from $\pi_a$ to $\pi_b$ in $G(\pi_1,\ldots,\pi_n)$, consecutive shadows intersect, whence
\[
\dist_*(U,V)\leq \sum_{i=1}^n\diam_*(\shadow(\pi_i))\leq 10\sum_{i=1}^n|\pi_i|\,.
\]
Consequently, $|w^{(i)}(\fX)|\leq e^{-s\,\dist_*(U,V)/10}\,\tilde w^{(i)}(\fX)$ for $\tilde w^{(i)}(\cdot):=|w^{(i)}(\cdot)|e^{s|\cdot|}$, which satisfies the hypotheses of \cref{eq:clusters-meeting-U} by \cref{eq:cmp-KP,eq:shadow-hyp}. Hence, by \cref{eq:clusters-meeting-U} for each of the four $\tilde w^{(i)}$,
\begin{align*}
\big|(\Sigma_1-\Sigma_2)-(\Sigma_3-\Sigma_4)\big|&=\bigg|\sum_{n\geq1}\frac1{n!}\sum_{\fX}\phi_\zeta(\fX)\big[w^{(1)}(\fX)-w^{(2)}(\fX)-w^{(3)}(\fX)+w^{(4)}(\fX)\big]\bigg|\\
&\leq4K|U|\,e^{-s\,\dist_*(U,V)/10}\,,
\end{align*}
yielding \cref{eq:log-X1-X2-log-X3-X4}, and thus completing the proof.
\end{proof}

\subsection{High temperature: proof of \texorpdfstring{\cref{prop:product}}{Proposition 3.1}}\label{subsec:proof-product}\mbox{}

\noindent\emph{Proof of Part~\labelcref{it:prod-KP}.} Since $\mu_\Lambda$ is a product measure and the supports of a compatible family are pairwise disjoint,
\[
\E_\Lambda\Big[\prod_{(\gamma,\eta)\in\cF}\one_{\{\sigma\restriction_\gamma=\eta\}}\Big]=\prod_{(\gamma,\eta)\in\cF}\mu_\Lambda(\sigma\restriction_\gamma=\eta)\,,
\]
so \cref{eq:polymer-rep} reads $Z_\Lambda^{\fg}=\sum_\cF\prod_{(\gamma,\eta)\in\cF}w(\gamma,\eta)$, the partition function of the polymer model defined in the proposition. Fix a polymer $(\gamma,\eta)$. If $(\gamma',\eta')$ is incompatible with it then $\dist(\gamma,\gamma')\leq1$ by \cref{eq:compatible-def}, so $\gamma'$ contains one of the at most $5|\gamma|$ sites at distance at most $1$ from $\gamma$; every $\gamma'\in\Gamma$ is connected, so that there are at most $(4e)^{n-1}$ possibilities for $\gamma'$ of a given size $n$ containing a given site (see the proof of \cref{cor:potential}), and $\sum_{\eta'\in\scS^{\gamma'}}\mu_\Lambda(\sigma\restriction_{\gamma'}=\eta')=1$. Hence, by \cref{eq:w-bound},
\begin{equation}\label{eq:prod-KP-sum}
\sum_{(\gamma',\eta')\,:\,\gamma'\not\perp\gamma} |w(\gamma',\eta')|\,e^{2|\gamma'|}\leq 5|\gamma|\sum_{n\geq1}(4e)^{n-1} e^{-(\Cst/10-2)n}\leq|\gamma|\,,
\end{equation}
provided that $(5+4e)e^{-(\Cst/10-2)}\leq1$, which holds for $\Cst>0$ large enough ($\Cst \geq48$ suffices). In particular, $w$ satisfies \cref{eq:textbook_bound} with $a(\gamma,\eta)=|\gamma|$. Since $w$ is real, the series in \cref{eq:cluster_expansion} thus converges absolutely to the logarithm of the positive quantity $Z_\Lambda^{\fg}$.

\noindent\emph{Proof of Part~\labelcref{it:prod-local}.} We apply \cref{lem:comparison} with $s=1$ to the family of polymers $(\gamma,\eta)$, $\gamma\in\Gamma_{\Lambda'}$, with $|(\gamma,\eta)|:=|\gamma|$ and $\shadow(\gamma,\eta):=\{y\in\Z^2:\dist_*(y,\gamma)\leq4|\gamma|\}$: as $\gamma$ is connected, $\diam_*(\shadow(\gamma,\eta))\leq9|\gamma|$, and incompatible polymers have intersecting shadows by \cref{eq:compatible-def}. Note that $\overline\gamma\subset\shadow(\gamma,\eta)$, as $\overline\gamma$ is contained in the smallest box containing $\gamma$. For a product measure on $\scS^{\Lambda'}$ or on $\scS^{\Lambda}$ (in the latter case the activities are extended by $0$ to the polymers of $\Lambda'$ not in $\Gamma_\Lambda$), the bound \cref{eq:prod-KP-sum} holds by the same proof, and is in fact \cref{eq:cmp-KP} with $s=1$ for the activities $w$; moreover, \cref{eq:shadow-hyp} holds with an absolute constant $K_0$: a polymer $(\gamma,\eta)$ with $|\gamma|=n$ and $\shadow(\gamma,\eta)\ni y$ has its support in $\{z:\dist_*(z,y)\leq4n\}$, so there are at most $(8n+1)^2(4e)^{n-1}$ such supports. Hence, $\sum_{\shadow(\pi)\ni y}|w(\pi)|e^{2|\pi|}\leq\sum_n(8n+1)^2(4e)^{n-1}e^{-(\Cst/10-2)n}\leq K_0$ for $\Cst$ large enough (the same $\Cst\geq 48$ from above suffices).

Fix $X$ and $\xi$ with $\mu_\Lambda(\sigma\restriction_{X}=\xi)>0$, and let $\scV^\xi$ be the specification with $\scV^\xi_{x,y}\equiv0$ and
\[
\scV^\xi_x(s):=\begin{cases}\scV_x(s) & \mbox{if $x\notin X$ or $s=\xi_x$}\,,\\
\infty & \mbox{if $x\in X$ and $s\neq\xi_x$}\,.\end{cases}
\]
Its Gibbs measure is $\mu_\Lambda^\xi=\mu_\Lambda(\cdot\given\sigma\restriction_X=\xi)$, and $Z_\Lambda(\scV^\xi)=Z_\Lambda\,\mu_\Lambda(\sigma\restriction_{X}=\xi)>0$. \Cref{thm:main} applied to $\scV^\xi$, together with Part~\labelcref{it:prod-KP}, gives
\begin{equation}\label{eq:prod-cond}
\E_\Lambda\big[\one_{\{\sigma\restriction_{X}=\xi\}}e^{\cH_\Lambda^\fg}\big]=\mu_\Lambda(\sigma\restriction_{X}=\xi)\,\E_\Lambda^\xi\big[e^{\cH_\Lambda^\fg}\big]=\mu_\Lambda(\sigma\restriction_{X}=\xi)\,\Xi_\Lambda^\xi\,,\qquad \Xi_\Lambda^\xi:=\sum_{\cF}\prod_{(\gamma,\eta)\in\cF}w^\xi(\gamma,\eta)\,,
\end{equation}
where $w^\xi(\gamma,\eta):=\sw^\xi(\gamma,\eta)\mu_\Lambda^\xi(\sigma\restriction_\gamma=\eta)$ and $\sw^\xi$ are the weights of \cref{thm:main} for $\scV^\xi$; therefore, $\mu_\Lambda^{\fg}(\sigma\restriction_{X}=\xi)=\mu_\Lambda(\sigma\restriction_{X}=\xi)\,\Xi_\Lambda^\xi/\Xi_\Lambda$, with $\Xi_\Lambda=Z_\Lambda^{\fg}$. Moreover $w^\xi\agreeoutside{X^c} w$: if $\shadow(\gamma,\eta)\cap X=\emptyset$ then $\overline\gamma\cap X=\emptyset$, so $\mu_\Lambda^\xi(\sigma\restriction_\gamma=\eta)=\mu_\Lambda(\sigma\restriction_\gamma=\eta)$, and $\sw^\xi(\gamma,\eta)=\sw(\gamma,\eta)$ since the latter depends on the specification only through the measure on $\Int(\gamma)$, which is the same in both cases. By Part~\labelcref{it:prod-KP}, applied to $\scV$ and to $\scV^\xi$, $\log\Xi_\Lambda$ and $\log\Xi_\Lambda^\xi$ are the sums of the corresponding series $\sum_n\fb_n$, which are the quantities $\Sigma_i$ compared in \cref{lem:comparison}.

Next, let $\Lambda'\supset\Lambda$ be finite and simply connected, and consider the four activities $w^\xi_{\Lambda'}$, $w_{\Lambda'}$, $w^\xi_\Lambda$, $w_\Lambda$ on the polymers of $\Lambda'$ (the last two extended by $0$). By the argument just given, $w^\xi_{\Lambda'}\agreeoutside{X^c} w_{\Lambda'}$ and $w^\xi_\Lambda\agreeoutside{X^c} w_\Lambda$. Moreover $w_{\Lambda'}$ and $w_\Lambda$ agree outside $V:=\Z^2\setminus\Lambda^\circ$, where $\Lambda^\circ:=\Lambda\setminus\partialvtxin\Lambda$: if $\shadow(\gamma,\eta)\subset\Lambda^\circ$ then $\gamma\in\Gamma_\Lambda$ (as $\gamma=\partialhat K$ with $K\subset\overline K=\overline\gamma\subset\Lambda$); the patterns contributing to $\sw(\gamma,\eta)$ are those with $B(x,r)\subset\overline\gamma$, the same collection for $\Lambda$ and for $\Lambda'$, carrying the same tilt since $\fg$ is defined in $\Z^2$ (or differs from such a tilt only on balls intersecting $\partialvtxin\Lambda$, as in \cref{rem:boundary-tilt}); and the measures on $\Int(\gamma)$ and on $\gamma$ are the same product measures. The same holds for $w^\xi_{\Lambda'}$ and $w^\xi_\Lambda$. Since $\mu_{\Lambda'}(\sigma\restriction_X=\xi)=\mu_\Lambda(\sigma\restriction_X=\xi)$ as $\scV$ is defined in $\Z^2$, the left-hand side of \cref{eq:prod-local-2} equals $\big|\log\frac{\Xi^\xi_{\Lambda'}}{\Xi_{\Lambda'}}-\log\frac{\Xi^\xi_{\Lambda}}{\Xi_{\Lambda}}\big|$ by \cref{eq:prod-cond}, which \cref{eq:log-X1-X2-log-X3-X4} with $U=X$, applied to the activities $w^\xi_{\Lambda'}$, $w_{\Lambda'}$, $w^\xi_\Lambda$, $w_\Lambda$ in this order, bounds by $4K_0|X|e^{-\dist_*(X,V)/10}$. As $\dist_*(X,V)\geq\dist_*(X,\Z^2\setminus\Lambda)-1\geq\frac12\dist(X,\Z^2\setminus\Lambda)-1$, this gives \cref{eq:prod-local-2} with $c:=1/20$ and $C'$ enlarged by an absolute factor. For \cref{eq:prod-local-3}, let $w^{\xi\xi'}$ be the activities of the specification with the spins frozen to $\xi$ on $X$ and to $\xi'$ on $Y$. By \cref{eq:prod-cond} for the three frozen specifications, and since $\mu_\Lambda(\sigma\restriction_X=\xi,\ \sigma\restriction_Y=\xi')=\mu_\Lambda(\sigma\restriction_X=\xi)\,\mu_\Lambda(\sigma\restriction_Y=\xi')$ for the product measure $\mu_\Lambda$,
\[
\log\frac{\mu_\Lambda^{\fg}(\sigma\restriction_{X}=\xi,\ \sigma\restriction_{Y}=\xi')}{\mu_\Lambda^{\fg}(\sigma\restriction_{X}=\xi)\,\mu_\Lambda^{\fg}(\sigma\restriction_{Y}=\xi')}=\log\Xi_\Lambda^{\xi\xi'}-\log\Xi_\Lambda^{\xi'}-\log\Xi_\Lambda^{\xi}+\log\Xi_\Lambda\,.
\]
By the same argument as above, $w^{\xi\xi'}\agreeoutside{X^c} w^{\xi'}$, $w^{\xi}\agreeoutside{X^c} w$, $w^{\xi\xi'}\agreeoutside{Y^c} w^{\xi}$ and $w^{\xi'}\agreeoutside{Y^c} w$, so \cref{eq:log-X1-X2-log-X3-X4} on $\Lambda$ with $U=X$ and $V=Y$, applied to the activities $w^{\xi\xi'}$, $w^{\xi'}$, $w^{\xi}$, $w$ in this order, bounds the right-hand side by $4K_0|X|e^{-\dist_*(X,Y)/10}$, giving \cref{eq:prod-local-3}.

Finally, \cref{eq:prod-local-2}, applied to two sets $\Lambda_n,\Lambda_m$ of a sequence $\Lambda_n\uparrow\Z^2$ with $\Lambda'$ a box containing both, shows that $\mu_{\Lambda_n}^{\fg}(\sigma\restriction_{X}=\xi)$ is Cauchy, with a limit that does not depend on the sequence (for $\xi$ with $\mu_\Lambda(\sigma\restriction_{X}=\xi)=0$ all these probabilities vanish), so the limits $\mu_{\infty}^{\fg}(\sigma\restriction_{X}=\xi)$ exist; being limits of consistent families of marginals, they define a probability measure $\mu_{\infty}^{\fg}$ on $\scS^{\Z^2}$, to which $\mu_\Lambda^{\fg}$ converges weakly. The bound \cref{eq:prod-local-3} passes to the limit. If $\scV$ and $\fg$ are translation invariant, then $\mu_{\Lambda+z}^{\fg}(\sigma\restriction_{X+z}=\xi(\cdot-z))=\mu_\Lambda^{\fg}(\sigma\restriction_{X}=\xi)$ for every $z\in\Z^2$, and the limit along the translated domains is the same, so $\mu_{\infty}^{\fg}$ is translation invariant. This completes the proof.
\qed 

\subsection{Low temperature: proof of \texorpdfstring{\cref{prop:lowT}}{Proposition 3.4}}\label{subsec:proof-lowT}

Recall the contour representation of the all-plus Ising model (see, e.g., \cite[\S3.7.2 and \S5.7.4]{FriedliVelenik18}). Let $E^*_\Lambda$ consist of the dual edges crossing primal edges with an endpoint in $\Lambda$, and put $\cX:=\{z\in\Z^2+(\frac12,\frac12):\dist_*(z,\Lambda)\leq\frac12\}$. Extending $\sigma$ by $+1$ outside $\Lambda$, its disagreement edges form a dual subgraph $\Gamma(\sigma)\subset E^*_\Lambda$, with $\cH_\Lambda^\tau(\sigma)=\beta|\Gamma(\sigma)|$. Its components, the \emph{contours}, are vertex-disjoint connected even subgraphs (every vertex has even degree). Conversely, such a family determines $\sigma_x=\prod_\Gamma(-1)^{N_\Gamma(x)}$, where $N_\Gamma(x)$ counts crossings of $\Gamma$ along a path from $x$ to infinity; its parity is path-independent. Connectedness of $\Z^2\setminus\Lambda$ ensures that the reconstructed spins are $+1$ outside $\Lambda$, giving a bijection. Write $\supp\Gamma$ for the vertex set of a contour; then $|\Gamma|\geq|\supp\Gamma|\geq4$.

\begin{remark}\label{rem:connected-complement}
Here connectedness of $\Z^2\setminus\Lambda$ is essential: if $z\notin\Lambda$ has all four neighbors in $\Lambda$, the four dual edges surrounding $z$ belong to $E^*_\Lambda$ but reconstruct $\sigma_z=-1$.
\end{remark}

For $A\subset\Lambda$, let $N(A)\subset\cX$ be the centers of unit squares with a corner in $A$; it is $*$-connected when $A$ is, and $|N(A)|\leq4|A|$. In either lattice there are at most $(8e)^{n-1}\leq64^n$ $*$-connected sets of size $n$ containing a given site, by the count in the proof of \cref{cor:potential} with $d=8$. Thus \cref{eq:effective-potential} and Cauchy--Schwarz give, for $\delta\geq5$,
\begin{equation}\label{eq:Psi-sqrt}
\|\Psi_A\|_\infty\leq e^{-\delta(|A|+1)}\,,\qquad
\sum_{A\ni y}\|\Psi_A\|_\infty^{1/2}\leq\Big(e^{-\delta}\sum_{n\geq1}64^ne^{-\delta n}\Big)^{1/2}\leq e^{-\delta/2}\,.
\end{equation}

\noindent \emph{Proof of Part~\labelcref{it:ising-KP}.}
Take $\beta_0=14$ and $C$ large enough that $\delta=\Cst/20-2\geq56$. The interior measures $\mu_{\Int(\gamma)}^\eta$ involve only the symmetric edge potentials, with boundary spins in $\gamma$ by \cref{eq:int-ext-dist}. Simultaneously flipping these spins and the interior spins, and using the symmetry of $\fg$, gives $\sw(\gamma,-\eta)=\sw(\gamma,\eta)$ by \cref{eq:Ec-def,eq:Ei-def,eq:w-def}, hence $\Psi_A(-\sigma)=\Psi_A(\sigma)$ by \cref{eq:Psi-def}.

Put $\vartheta_A:=e^{\Psi_A-\Psi_A(\underline{1})}-1$, so $\|\vartheta_A\|_\infty\leq4\|\Psi_A\|_\infty$. Since $A$ is connected, the disagreement edges along paths in $A$ determine $\sigma\restriction_A$ up to a global flip. Their dual edges have endpoints in $N(A)$; hence $\vartheta_A$ depends only on contours meeting $N(A)$, and vanishes when none do. Expanding \cref{eq:effective-potential} therefore gives
\[
Z_\Lambda^\tau Z_\Lambda^{\fg}
=e^{\sum_A\Psi_A(\underline{1})}\sum_{\{\Gamma_j\}}\sum_{\cA}\prod_j e^{-\beta|\Gamma_j|}\prod_{A\in\cA}\vartheta_A(\sigma)\,,
\]
where the contour family determines $\sigma$ and $\cA$ ranges over collections of connected subsets of $\Lambda$. Call the contours and the sets in $\cA$ \emph{atoms}, with supports $\supp\Gamma$ and $N(A)$, respectively. Group atoms into components $P$ of their support-intersection graph. Each component has a $*$-connected support $S$, distinct components have disjoint supports, and every factor $\vartheta_A$ depends only on contours in its component. Thus, writing $\sigma_P$ for the configuration reconstructed from the contours in $P$,
\[
F_P:=\prod_{\Gamma\in P}e^{-\beta|\Gamma|}\prod_{A\in P}\vartheta_A(\sigma_P)\,,\qquad
W(S):=\sum_{P\,:\,\supp P=S}F_P\,,\qquad
Z_\Lambda^\tau Z_\Lambda^{\fg}=e^{\sum_A\Psi_A(\underline{1})}\sum_{\{S_i\}}\prod_iW(S_i)\,.
\]
Here $P$ ranges over sets of atoms with connected intersection graph and vertex-disjoint contours, and the last sum is over families of pairwise disjoint supports. The grouping is reversible, so this is the required polymer representation.

Set $w_\Gamma:=e^{-\beta|\Gamma|}$ and $w_A:=4\|\Psi_A\|_\infty$, so $|F_P|\leq\prod_{\alpha\in P}w_\alpha$. We have
\begin{equation}\label{eq:atom-hyp}
w_\alpha\leq e^{-14|\supp\alpha|}\,,\qquad
\sum_{\alpha\,:\,\supp\alpha\ni z}w_\alpha^{1/2}\leq1\,.
\end{equation}
The first bound follows from $|\Gamma|\geq|\supp\Gamma|$, $|N(A)|\leq4|A|$ and \cref{eq:Psi-sqrt}. For the second, a contour of $n$ edges through $z$ has an Euler circuit starting at $z$, so there are at most $4^n$ such contours; and $N(A)\ni z$ forces $A$ to contain a corner of the square centered at $z$. Hence,
\[
\sum_{\alpha\,:\,\supp\alpha\ni z}w_\alpha^{1/2}
\leq\sum_{n\geq4}(4e^{-\beta/2})^n+8e^{-\delta/2}\leq1\,.
\]
Splitting each $w_\alpha$ into two square roots and then dropping all constraints on $P$ except $\supp\alpha\subset S$, we obtain
\begin{equation}\label{eq:W-bound}
|W(S)|\leq e^{-7|S|}\prod_{\alpha\,:\,\supp\alpha\subset S}(1+w_\alpha^{1/2})
\leq e^{-7|S|}\exp\Big(\sum_{z\in S}\sum_{\alpha\,:\,\supp\alpha\ni z}w_\alpha^{1/2}\Big)\leq e^{-6|S|}\,.
\end{equation}
Consequently,
\begin{equation}\label{eq:W-KP}
\sum_{S'\,:\,S'\cap S\neq\emptyset}|W(S')|e^{(1+1/10)|S'|}
\leq|S|\sum_{n\geq1}64^ne^{-4.9n}\leq|S|\,.
\end{equation}
This is \cref{eq:cmp-KP} with $s=1/10$; since $W$ is real, the series $\sum_n\fb_n$ for the activities $W$ is the logarithm of the positive quantity $e^{-\sum_A\Psi_A(\underline{1})}Z_\Lambda^\tau Z_\Lambda^{\fg}$, proving Part~\labelcref{it:ising-KP}.

\noindent \emph{Proof of Part~\labelcref{it:ising-local}.}
For $T\subset\Lambda$, write $\sigma_T:=\prod_{x\in T}\sigma_x$ as well as $\epsilon_T(\Gamma):=(-1)^{\sum_{x\in T}N_\Gamma(x)}$. Since $\sigma_T=\prod_\Gamma\epsilon_T(\Gamma)$, inserting it in the preceding expansion gives
\begin{equation}\label{eq:ising-signed}
\Xi_{T,\Lambda}:=\sum_\sigma\sigma_T e^{-\beta|\Gamma(\sigma)|}\prod_A(1+\vartheta_A(\sigma))
=\sum_{\{S_i\}}\prod_iW_T(S_i)\,,\qquad
W_T(S):=\sum_{P\,:\,\supp P=S}F_P\prod_{\Gamma\in P}\epsilon_T(\Gamma)\,.
\end{equation}
Notice that \cref{eq:W-bound,eq:W-KP} also hold for $W_T$. Now let $\shadow(S):=\{y\in\Z^2:\dist_*(y,Q_S)\leq4|S|\}$, where $Q_S$ consists of the corners of squares centered at sites of $S$. These have diameter at most $9|S|$, and incompatible polymers have intersecting shadows. A contour enclosing $x$ has $x$ in its bounding box, so changing $T$ on $U$ changes $W_T$ only on polymers whose shadows meet $U$. Moreover, for $\Lambda'\supset\Lambda$ and globally defined $\fg$, the activities $W_{T,\Lambda'}$ and $W_{T,\Lambda}$ agree outside $V:=\Z^2\setminus\Lambda^\circ$, where $\Lambda^\circ:=\Lambda\setminus\partialvtxin\Lambda$, extending the latter by zero. Indeed, let $\shadow(S)\subset\Lambda^\circ$. The atoms with support in $S$ are then the same in both domains: a contour with vertex set in $S$ consists of dual edges crossing primal edges with endpoints in $Q_S\subset\Lambda^\circ$, which belong to $E^*_\Lambda$, and a connected $A$ with $N(A)\subset S$ satisfies $A\subset Q_S$. For such $A$, every tuple in \cref{eq:Psi-def} has $\gamma_i\subset A$, so $\overline{\gamma_i}$ lies in the bounding box of $A$, within $\ell_\infty$-distance $\diam_*(A)<|A|\leq|N(A)|\leq|S|$ of $Q_S$, hence in $\shadow(S)$, and $\gamma_i\in\Gamma_\Lambda$ as in the proof of \cref{prop:product}\labelcref{it:prod-local}. The weight $\sw(\gamma_i,\eta)$ depends only on the tilt in $\overline{\gamma_i}$ and on the Ising measure in $\Int(\gamma_i)$ with boundary spins $\eta$ on $\gamma_i$, by \cref{eq:int-ext-dist}, so the all-plus boundary condition does not enter, and $\Psi_A$ is the same for $\Lambda$ and $\Lambda'$. Finally, a size-$n$ polymer whose shadow contains $y$ has a site within distance $4n+1$ of $y$, so \cref{eq:W-bound} and the bound $64^n$ imply \cref{eq:shadow-hyp} with an absolute $K_0$ and $s=1/10$.

Let $\widetilde\mu_\Lambda\propto e^{\sum_A\Psi_A}\mu_\Lambda$. The convergent real cluster expansion gives $\Xi_{T,\Lambda}>0$, while \cref{eq:ising-signed} gives
\[
0<M_{T,\Lambda}:=\frac{\Xi_{T,\Lambda}}{\Xi_{\emptyset,\Lambda}}=\widetilde\mu_\Lambda(\sigma_T)\leq1\,.
\]
Since $|a-b|\leq|\log a-\log b|$ for $a,b\in(0,1]$, \cref{lem:comparison} applied to $W_{T,\Lambda'},W_{\Lambda'},W_{T,\Lambda},W_\Lambda$, and to $W_{T\cup T'},W_{T'},W_T,W$ for disjoint $T,T'$, yields
\begin{equation}\label{eq:ising-moments}
\begin{split}
|M_{T,\Lambda'}-M_{T,\Lambda}|&\leq K|T|e^{-c\,\dist_*(T,\Z^2\setminus\Lambda)}\,,\\
|M_{T\cup T',\Lambda}-M_{T,\Lambda}M_{T',\Lambda}|&\leq K|T|e^{-c\,\dist_*(T,T')}\,.
\end{split}
\end{equation}
Here and below $K,c>0$ are absolute constants, enlarged or decreased as needed. For a function $f=\sum_T\hat f(T)\sigma_T$, write $\|f\|_{\mathrm F}:=\sum_T|\hat f(T)|e^{|T|}$. Summing \cref{eq:ising-moments} and using $|T|\leq e^{|T|}$ gives, for functions supported on $R,R'$ respectively,
\begin{equation}\label{eq:ising-functions}
\begin{split}
|\widetilde\mu_{\Lambda'}(f)-\widetilde\mu_\Lambda(f)|&\leq K\|f\|_{\mathrm F}e^{-c\,\dist_*(R,\Z^2\setminus\Lambda)}\,,\\
|\widetilde\mu_\Lambda(ff')-\widetilde\mu_\Lambda(f)\widetilde\mu_\Lambda(f')|&\leq K\|f\|_{\mathrm F}\|f'\|_{\mathrm F}e^{-c\,\dist_*(R,R')}
\end{split}
\end{equation}
(the second bound is for disjoint $R,R'$).

To pass to $\mu_\Lambda^{\fg}$, freeze the spins on $X$ to $\xi$ and denote the resulting effective potential by $\Psi_A^\xi$. Put $\theta_A:=e^{\Psi_A^\xi-\Psi_A}-1$, calling $A$ \texttt{bad} when $\theta_A\not\equiv0$. By \cref{eq:Psi-def}, whose weights depend on the specification only within the hulls $\overline{\gamma_i}$, which lie within $\ell_\infty$-distance $|A|$ of $A$, a \texttt{bad} $A$ has a site within $\ell_\infty$-distance $|A|$ of $X$, so that $A$ and these hulls lie within $\ell_\infty$-distance $2|A|$ and $3|A|$ of $X$, respectively; moreover, \cref{eq:effective-potential} gives $\|\theta_A\|_\infty\leq4e^{-\delta(|A|+1)}$. Applying the partition-function identity to the frozen specification gives exactly
\begin{equation}\label{eq:ising-cond}
\mu_\Lambda^{\fg}(\sigma\restriction_X=\xi)=\widetilde\mu_\Lambda(h_\xi)\,,\qquad
h_\xi:=\one_{\{\sigma\restriction_X=\xi\}}\prod_{A\texttt{ bad}}(1+\theta_A)\,.
\end{equation}
For $d\geq20$, let $h_\xi^{(d)}$ retain only factors with $|A|\leq d/20$. There are at most $|X|(8n+1)^264^n$ \texttt{bad} sets of size $n$, since each contains a site within distance $4n$ of $X$. Hence, we can infer that
\begin{equation}\label{eq:ising-trunc}
\|h_\xi-h_\xi^{(d)}\|_\infty\leq e^{K|X|}e^{-cd}\,,\qquad
\|h_\xi^{(d)}\|_{\mathrm F}\leq e^{K|X|}\,.
\end{equation}
Indeed, the sum of $\|\theta_A\|_\infty$ over \texttt{bad} sets is at most $K|X|$, and its tail above $d/20$ is at most $K|X|e^{-cd}$, proving the first bound by telescoping the product. For the second bound in \cref{eq:ising-trunc}, the Fourier norm is submultiplicative, $\|\one_{\{\sigma\restriction_X=\xi\}}\|_{\mathrm F}\leq e^{|X|}$, and Cauchy--Schwarz and Parseval give $\|\theta_A\|_{\mathrm F}\leq e^{|A|}2^{|A|/2}\|\theta_A\|_\infty$; the same count bounds the product of $1+\|\theta_A\|_{\mathrm F}$ by $e^{K|X|}$.

Recall that a retained \texttt{bad} set $A$, and the hulls through which $\Psi_A$ and $\Psi_A^\xi$ depend on the specification, lie within $\ell_\infty$-distance $3|A|\leq3d/20$ of $X$. For \cref{eq:ising-local-2}, take $d=\dist_*(X,\Z^2\setminus\Lambda)$. These hulls then lie in $\Lambda^\circ$, so, by the locality of $\Psi_A$ established above, applied to the Ising specification and to the frozen one, the correction $\theta_A$ of every retained \texttt{bad} set is the same in $\Lambda$ and $\Lambda'$. Thus $h_\xi^{(d)}$ is common to both domains and supported at $\ell_\infty$-distance at least $d/2$ from $\Z^2\setminus\Lambda$, and the first bound in \cref{eq:ising-functions}, together with \cref{eq:ising-trunc} in both domains, gives
\[
\big|\mu_{\Lambda'}^{\fg}(\sigma\restriction_X=\xi)-\mu_\Lambda^{\fg}(\sigma\restriction_X=\xi)\big|\leq\big|\widetilde\mu_{\Lambda'}(h_\xi^{(d)})-\widetilde\mu_\Lambda(h_\xi^{(d)})\big|+2e^{K|X|}e^{-cd}\leq Ke^{K|X|}e^{-cd/2}+2e^{K|X|}e^{-cd}\,.
\]
For \cref{eq:ising-local-3}, take $d=\dist_*(X,Y)$ and truncate the corrections for $\xi$, $\xi'$ and $(\xi,\xi')$. A retained \texttt{bad} set for $(\xi,\xi')$ lies within $\ell_\infty$-distance $3d/20$ of $X$ or of $Y$, and freezing the spins on the other set, which is at $\ell_\infty$-distance at least $d-3d/20$ from the hulls through which its correction depends on the specification, does not affect that correction; thus $h_{\xi\xi'}^{(d)}=h_\xi^{(d)}h_{\xi'}^{(d)}$, with the two factors supported within $\ell_\infty$-distance $3d/20$ of $X$ and of $Y$, respectively, hence at $\ell_\infty$-distance at least $d/2$ from one another. Since $\widetilde\mu_\Lambda(h_\xi)=\mu_\Lambda^\fg(\sigma\restriction_X=\xi)\leq1$ by \cref{eq:ising-cond}, replacing $h_\xi$, $h_{\xi'}$ and $h_{\xi\xi'}$ by their truncations changes the left-hand side of \cref{eq:ising-local-3} by at most $3e^{K(|X|+|Y|)}e^{-cd}$ by \cref{eq:ising-trunc}, after which the second bound in \cref{eq:ising-functions} bounds it by $Ke^{K(|X|+|Y|)}e^{-cd/2}$. In both cases $d<20$ is trivial after increasing $C'$, and $\dist_*\geq\frac12\dist$ converts to the stated distances. Weak convergence and translation invariance follow from \cref{eq:ising-local-2} exactly as in \cref{prop:product}\labelcref{it:prod-local}.
\qed

\appendix
\section{Proofs of \texorpdfstring{\cref{lem:markov}}{the domain Markov lemma} and \texorpdfstring{\cref{clm:ursell-sign}}{the Ursell functions}}\label{app:markov}
\begin{proof}[Proof of \cref{lem:markov}]
Let $\sigma\in\cE_\xi$, and split the summands in $\cH_V^\tau(\sigma)$ per \cref{eq:gibbs-def} according to whether they involve $\sigma\restriction_U$. The edges of $E(V)$ with both endpoints in $V\setminus U$, the sites of $V\setminus U$, and the edges of~$\partialedge V$ whose endpoint in $V$ belongs to $V\setminus U$ contribute terms depending only on $\xi$ (recall that $\xi\restriction_{V^c}=\tau$); we denote their sum by $A_\xi\in\R\cup\{+\infty\}$, which does not depend on $\sigma$. The remaining terms are indexed by the edges of $E(U)$, the sites of $U$, and the edges joining $U$ to $V\setminus U$ or to $V^c$; the last two families together form $\partialedge U$, and on each such edge $xz$ with $x\in U$ the term is $\scV_{x,z}(\sigma_x,\xi_z)$. Comparing with \cref{eq:gibbs-def} for $U$, we find that
\[
\cH_V^\tau(\sigma) = A_\xi + \cH_U^\xi(\sigma\restriction_U)\qquad\mbox{for all $\sigma\in\cE_\xi$}\,;
\]
thus,
summing $e^{-\cH_V^\tau(\sigma)}$ over $\sigma\in\cE_\xi$ gives $\mu_V^\tau(\cE_\xi)=e^{-A_\xi}Z_U^\xi/Z_V^\tau$, which is positive by assumption. In particular,
\[
\mu_V^\tau(\cE_\xi)\,\E_U^\xi[F]
=\frac{e^{-A_\xi}}{Z_V^\tau}\sum_{\tilde\sigma\in\scS^U}F(\tilde\sigma)e^{-\cH_U^\xi(\tilde\sigma)}=
\frac{1}{Z_V^\tau} \sum_{\sigma\in\cE_\xi} F(\sigma\restriction_U) e^{-\cH_V^\tau(\sigma)} = 
\E_V^\tau\big[F(\sigma\restriction_U)\,\one_{\cE_\xi}\big]  \,,
\]
as claimed.

For the last assertion, no edge of $\Z^2$ joins $U_j$ to $U_k$ for $j\neq k$, so $E(U)$ is the disjoint union of $E(U_1),\ldots,E(U_m)$, while $\partialedge U$ is the disjoint union of $\partialedge U_1,\ldots,\partialedge U_m$. In particular $\partialvtxout U_k\subset U^c$, so $\cH_{U_k}^\xi$, $Z_{U_k}^\xi$ and $\mu_{U_k}^\xi$ are defined via the convention following \cref{eq:gibbs-def}, and \cref{eq:gibbs-def} for $U$ reads $\cH_U^\xi(\sigma)=\sum_{k=1}^m \cH_{U_k}^\xi(\sigma\restriction_{U_k})$ for all $\sigma\in\scS^U$. Thus $Z_U^\xi=\prod_{k=1}^m Z_{U_k}^\xi$, whence every $Z_{U_k}^\xi$ is positive, and $\mu_U^\xi=\bigotimes_{k=1}^m\mu_{U_k}^\xi$, which concludes the proof.
\end{proof}

\begin{proof}[Proof of \cref{clm:ursell-sign}]
For a finite multigraph $G$ (possibly with loops and multiple edges), let
\[
t(G):=(-1)^{|V(G)|-1}\sum_{\substack{H\subset G\\\mbox{conn., span.}}}(-1)^{|E(H)|}\,,
\]
the sum running over all subgraphs $H$ of $G$ that are spanning (i.e., $V(H)=V(G)$) and connected. By its definition in \cref{eq:ursell-def}, in the special case $\zeta(\pi_1,\pi_2) = \one_{\{\pi_1\perp \pi_2\}}$ we have that
\[ \phi_\zeta(\pi_1,\ldots,\pi_n) = (-1)^{n-1} t(G)\,,
\] 
with the \emph{incompatibility graph} $G$ whose vertices are $\{v_1,\ldots,v_n\}$ and edges are $\{v_i v_j\,:\; \pi_i \not\perp \pi_j\}$. 
The claim is then that $t(G)\geq 0$, where one may assume that $G$ is connected (or else $t(G)=0$), and that $t(K_n)=(n-1)!$, where $K_n$ is the complete (simple) graph (recalling that $\pi \not\perp \pi$).

In terms of the Tutte polynomial (see \cite[Ch.~13]{Biggs_AGT}), summing over spanning subgraphs $H$ of $G$,
\[ T_G(x,y)=\sum_{H\subset G}(x-1)^{k(H)-k(G)}(y-1)^{|E(H)|-|V(G)|+k(H)}\,,\] where $k(\cdot)$ is the number of connected components, one has $t(G)=T_G(1,0)$ for a connected $G$, since at $x=1$ only the spanning connected subgraphs survive. The two sought facts are classical:
\begin{enumerate}[(a)]
\item $T_G(x,y)$ has nonnegative coefficients by Tutte's spanning-tree expansion \cite[\S3]{Tutte54} (see the definition of $T_G$ in \cite[Ch.~13]{Biggs_AGT} and equivalence in \cite[Thm.~13.9]{Biggs_AGT}); in particular, $T_G(1,0)\geq 0$.
\item $T_G(1,0)$ is the number of acyclic orientations of $G$ with a prescribed vertex $v$ as a unique source \cite[Thm.~7.3]{GZ83}. For $G=K_n$, these are the $(n-1)!$ linear orderings starting from $v$.
\end{enumerate}

For completeness, we also provide the simple direct proofs of these two facts.

We show that $t(G)\geq 0$ for any multigraph $G$ by induction on $|E(G)|$, starting from the trivial base case $E(G)=\emptyset$ wherein $t(G)=1$ if $|V(G)|=1$ and $t(G)=0$ if $|V(G)|>1$. 
If $G$ contains some loop $e$ then immediately $t(G)=0$: the bijection from $H\ni e$ to $H':=H\setminus\{e\}$ has $(-1)^{|E(H')|}=(-1)^{|E(H)|+1}$, so they cancel out in $t(G)$.
For a loopless $G$, let $e=v_i v_j\in E(G)$, and write $t(G) = t_e(G) + t_e^{\neg}(G)$ where
\[ t_e(G) = (-1)^{|V(G)|-1}\sum_{H \ni e} (-1)^{|E(H)|}\,,\quad
t_e^{\neg}(G) = (-1)^{|V(G)|-1}\sum_{H \not\ni e} (-1)^{|E(H)|}\,.\]
We argue that $t_e(G)\geq 0$ and $t_e^\neg(G)\geq 0$. 
For the latter, observe that $t_e^\neg(G)$ is nothing but $t(G')$ for the graph $G':= G \setminus\{e\}$, and hence nonnegative by the induction hypothesis. For the former, consider $G''$ obtained by contracting $e$ (deleting it and then merging the vertices $v_i,v_j$ to a single one). Contracting $e$ is also a bijection from connected spanning subgraphs $H\ni e$ of $G$ to connected spanning subgraphs $H''$ of $G''$, and $(-1)^{|V(G)|+|E(H)|}=(-1)^{|V(G'')|+|E(H'')|}$ as it deletes one vertex and one edge. Thus, $t_e(G) = t(G'')$, which again is nonnegative by the induction hypothesis.

To show that $t(K_n)=(n-1)!$, we argue by induction on $n$, starting from the trivial case $n=1$. For $n\geq2$, group the spanning connected subgraphs $H$ of $K_n$ according to the set $N_1\neq\emptyset$ of neighbors of $v_1$ in $H$. For a given $N_1$, its internal edge configuration in $H$ does not affect the connectivity of $H$, as $N_1\cup\{v_1\}$ is connected through $v_1$; summing over it thus produces the factor $\sum_{j}\binom{m}{j}(-1)^j=\one_{\{m=0\}}$ for $m=\binom{|N_1|}2$, so that only $|N_1|=1$ contributes. If $N_1=\{w\}$, then $H$ is connected if and only if $H-v_1$ is a spanning connected subgraph of $K_{n-1}$, and $H$ has one more vertex and one more edge than $H-v_1$, so the signs in the definition of $t(G)$ cancel. Summing over the $n-1$ choices of $w$ and using the induction hypothesis, $t(K_n)=(n-1)\,t(K_{n-1})=(n-1)!$.
\end{proof}

\subsection*{Acknowledgments}
The research of BL was supported by the grant \textsc{anr-18-ce40-0033 dimers}. The research of EL was supported by the NSF grant \textsc{dms}-2451083.

The heart of the argument behind the main theorem was key to the renormalization group analysis used by the authors in \cite[Prop.~8.11]{LaslierLubetzky24}. From the referee's feedback on the nontrivial aspect of having, in that paper's system of Grassmann variables, coefficients that decay with the radius of the associated pattern as opposed to its area, we decided to develop this idea in the more general setting of this paper, while also dispensing with the restriction on the sign of the interaction.

AI assisted us in preparing this manuscript in the following way. We gave \textsc{Claude Fable~5.1} a file containing the proofs of the main results, and used it to audit those and polish their exposition (it found, and readily fixed, various small issues in the combinatorics associated to our polymers, mostly relating to adjacency vs.\ $*$-adjacency conventions affecting their connectedness). We further used it to review the related literature, and to assist with optimizing the proofs, and constants, of the cluster expansion applications that appear in \cref{prop:product,prop:lowT}.

\bibliographystyle{abbrv}
\bibliography{ce_longrange}

@article {GMT20,
    AUTHOR = {Giuliani, Alessandro and Mastropietro, Vieri and Toninelli,
              Fabio Lucio},
     TITLE = {Non-integrable dimers: universal fluctuations of tilted height
              profiles},
   JOURNAL = {Comm. Math. Phys.},
  FJOURNAL = {Communications in Mathematical Physics},
    VOLUME = {377},
      YEAR = {2020},
    NUMBER = {3},
     PAGES = {1883--1959},
      ISSN = {0010-3616},
   MRCLASS = {82B20 (37K60 82B28)},
  MRNUMBER = {4121614},
MRREVIEWER = {Arvind Ayyer},
       DOI = {10.1007/s00220-020-03760-x},
       URL = {https://doi.org/10.1007/s00220-020-03760-x},
}

@article {GMT17,
    AUTHOR = {Giuliani, Alessandro and Mastropietro, Vieri and Toninelli,
              Fabio Lucio},
     TITLE = {Height fluctuations in interacting dimers},
   JOURNAL = {Ann. Inst. Henri Poincar\'{e} Probab. Stat.},
  FJOURNAL = {Annales de l'Institut Henri Poincar\'{e} Probabilit\'{e}s et
              Statistiques},
    VOLUME = {53},
      YEAR = {2017},
    NUMBER = {1},
     PAGES = {98--168},
      ISSN = {0246-0203},
   MRCLASS = {82B28 (82B20)},
  MRNUMBER = {3606736},
       DOI = {10.1214/15-AIHP710},
       URL = {https://doi.org/10.1214/15-AIHP710},
}

@book{FriedliVelenik18,
    AUTHOR = {Friedli, S. and Velenik, Y.},
     TITLE = {Statistical mechanics of lattice systems},
      NOTE = {A concrete mathematical introduction},
 PUBLISHER = {Cambridge University Press, Cambridge},
      YEAR = {2018},
     PAGES = {xix+622},
      ISBN = {978-1-107-18482-4},
   MRCLASS = {82-01 (82B05)},
  MRNUMBER = {3752129},
}

@article{KoteckyPreiss86,
    AUTHOR = {Koteck\'{y}, R. and Preiss, D.},
     TITLE = {Cluster expansion for abstract polymer models},
   JOURNAL = {Comm. Math. Phys.},
  FJOURNAL = {Communications in Mathematical Physics},
    VOLUME = {103},
      YEAR = {1986},
    NUMBER = {3},
     PAGES = {491--498},
       DOI = {10.1007/BF01211762},
       URL = {https://doi.org/10.1007/BF01211762},
}

@article{FernandezProcacci07,
    AUTHOR = {Fern\'{a}ndez, Roberto and Procacci, Aldo},
     TITLE = {Cluster expansion for abstract polymer models. {N}ew bounds from an old approach},
   JOURNAL = {Comm. Math. Phys.},
  FJOURNAL = {Communications in Mathematical Physics},
    VOLUME = {274},
      YEAR = {2007},
    NUMBER = {1},
     PAGES = {123--140},
       DOI = {10.1007/s00220-007-0279-2},
       URL = {https://doi.org/10.1007/s00220-007-0279-2},
}

@article{Timar13,
    AUTHOR = {Tim\'{a}r, \'{A}d\'{a}m},
     TITLE = {Boundary-connectivity via graph theory},
   JOURNAL = {Proc. Amer. Math. Soc.},
  FJOURNAL = {Proceedings of the American Mathematical Society},
    VOLUME = {141},
      YEAR = {2013},
    NUMBER = {2},
     PAGES = {475--480},
       DOI = {10.1090/S0002-9939-2012-11333-4},
       URL = {https://doi.org/10.1090/S0002-9939-2012-11333-4},
}

@article{Hammond06,
    AUTHOR = {Hammond, Alan},
     TITLE = {Greedy lattice animals: geometry and criticality},
   JOURNAL = {Ann. Probab.},
  FJOURNAL = {The Annals of Probability},
    VOLUME = {34},
      YEAR = {2006},
    NUMBER = {2},
     PAGES = {593--637},
      ISSN = {0091-1798,2168-894X},
   MRCLASS = {60K35 (82B43)},
  MRNUMBER = {2223953},
MRREVIEWER = {Rongfeng\ Sun},
       DOI = {10.1214/009117905000000693},
       URL = {https://doi.org/10.1214/009117905000000693},
}

@incollection{Kesten86,
    AUTHOR = {Kesten, Harry},
     TITLE = {Aspects of first passage percolation},
 BOOKTITLE = {\'{E}cole d'\'{e}t\'{e} de probabilit\'{e}s de {S}aint-{F}lour, {XIV}---1984},
    SERIES = {Lecture Notes in Math.},
    VOLUME = {1180},
     PAGES = {125--264},
 PUBLISHER = {Springer, Berlin},
      YEAR = {1986},
       DOI = {10.1007/BFb0074919},
}

@misc{LaslierLubetzky24,
    AUTHOR = {Laslier, Beno\^{\i}t and Lubetzky, Eyal},
     TITLE = {Tilted {S}olid-{O}n-{S}olid is liquid: scaling limit of {SOS} with a potential on a slope},
      NOTE = {Preprint, \texttt{arXiv:2409.08745}},
      YEAR = {2024},
}

@article{Zahradnik84,
    AUTHOR = {Zahradn\'{\i}k, Milo\v{s}},
     TITLE = {An alternate version of {P}irogov-{S}inai theory},
   JOURNAL = {Comm. Math. Phys.},
  FJOURNAL = {Communications in Mathematical Physics},
    VOLUME = {93},
      YEAR = {1984},
    NUMBER = {4},
     PAGES = {559--581},
       DOI = {10.1007/BF01212295},
       URL = {https://doi.org/10.1007/BF01212295},
}

@incollection{Dobrushin96,
    AUTHOR = {Dobrushin, R. L.},
     TITLE = {Estimates of semi-invariants for the {I}sing model at low temperatures},
 BOOKTITLE = {Topics in statistical and theoretical physics},
    SERIES = {Amer. Math. Soc. Transl. Ser. 2},
    VOLUME = {177},
     PAGES = {59--81},
 PUBLISHER = {Amer. Math. Soc., Providence, RI},
      YEAR = {1996},
}

@article{ABEH25,
    AUTHOR = {Affonso, Lucas and Bissacot, Rodrigo and Endo, Eric O. and Handa, Satoshi},
     TITLE = {Long-range {I}sing models: contours, phase transitions and decaying fields},
   JOURNAL = {J. Eur. Math. Soc. (JEMS)},
  FJOURNAL = {Journal of the European Mathematical Society (JEMS)},
    VOLUME = {27},
      YEAR = {2025},
    NUMBER = {4},
     PAGES = {1679--1714},
       DOI = {10.4171/JEMS/1529},
       URL = {https://doi.org/10.4171/JEMS/1529},
      NOTE = {arXiv:2105.06103},
}

@article{GruberKunz71,
    AUTHOR = {Gruber, C. and Kunz, H.},
     TITLE = {General properties of polymer systems},
   JOURNAL = {Comm. Math. Phys.},
  FJOURNAL = {Communications in Mathematical Physics},
    VOLUME = {22},
      YEAR = {1971},
    NUMBER = {2},
     PAGES = {133--161},
       DOI = {10.1007/BF01651334},
       URL = {https://doi.org/10.1007/BF01651334},
}

@article{Mayer37,
    AUTHOR = {Mayer, Joseph E.},
     TITLE = {The statistical mechanics of condensing systems. {I}},
   JOURNAL = {J. Chem. Phys.},
  FJOURNAL = {The Journal of Chemical Physics},
    VOLUME = {5},
      YEAR = {1937},
    NUMBER = {1},
     PAGES = {67--73},
       DOI = {10.1063/1.1749933},
       URL = {https://doi.org/10.1063/1.1749933},
}

@article{Ursell27,
    AUTHOR = {Ursell, H. D.},
     TITLE = {The evaluation of {G}ibbs' phase-integral for imperfect gases},
   JOURNAL = {Proc. Cambridge Philos. Soc.},
  FJOURNAL = {Proceedings of the Cambridge Philosophical Society},
    VOLUME = {23},
      YEAR = {1927},
    NUMBER = {6},
     PAGES = {685--697},
       DOI = {10.1017/S0305004100011191},
       URL = {https://doi.org/10.1017/S0305004100011191},
}

@article{MinlosSinai67,
    AUTHOR = {Minlos, R. A. and Sina\u{\i}, Ya. G.},
     TITLE = {The phenomenon of ``separation of phases'' at low temperatures
              in certain lattice models of a gas. {I}},
   JOURNAL = {Math. USSR-Sb.},
  FJOURNAL = {Mathematics of the USSR-Sbornik},
    VOLUME = {2},
      YEAR = {1967},
    NUMBER = {3},
     PAGES = {335--395},
  MRNUMBER = {213133},
       DOI = {10.1070/SM1967v002n03ABEH002345},
       URL = {https://doi.org/10.1070/SM1967v002n03ABEH002345},
      NOTE = {Translated from Mat. Sb. (N.S.) \textbf{73(115)} (1967), no.~3, 375--448},
}

@article{Dobrushin72,
    AUTHOR = {Dobrushin, R. L.},
     TITLE = {Gibbs state describing coexistence of phases for a
              three-dimensional {I}sing model},
   JOURNAL = {Theory Probab. Appl.},
  FJOURNAL = {Theory of Probability and its Applications},
    VOLUME = {17},
      YEAR = {1972},
    NUMBER = {4},
     PAGES = {582--600},
  MRNUMBER = {421546},
       DOI = {10.1137/1117073},
       URL = {https://doi.org/10.1137/1117073},
      NOTE = {Translated from Teor. Verojatnost. i Primenen. \textbf{17} (1972), no.~4, 619--639},
}

@article{PirogovSinai75,
    AUTHOR = {Pirogov, S. A. and Sina\u{\i}, Ya. G.},
     TITLE = {Phase diagrams of classical lattice systems},
   JOURNAL = {Theoret. and Math. Phys.},
  FJOURNAL = {Theoretical and Mathematical Physics},
    VOLUME = {25},
      YEAR = {1975},
    NUMBER = {3},
     PAGES = {1185--1192},
  MRNUMBER = {676316},
       DOI = {10.1007/BF01040127},
       URL = {https://doi.org/10.1007/BF01040127},
      NOTE = {Translated from Teoret. Mat. Fiz. \textbf{25} (1975), no.~3, 358--369},
}

@article{PirogovSinai76,
    AUTHOR = {Pirogov, S. A. and Sina\u{\i}, Ya. G.},
     TITLE = {Phase diagrams of classical lattice systems. {C}ontinuation},
   JOURNAL = {Theoret. and Math. Phys.},
  FJOURNAL = {Theoretical and Mathematical Physics},
    VOLUME = {26},
      YEAR = {1976},
    NUMBER = {1},
     PAGES = {39--49},
  MRNUMBER = {676499},
       DOI = {10.1007/BF01038255},
       URL = {https://doi.org/10.1007/BF01038255},
      NOTE = {Translated from Teoret. Mat. Fiz. \textbf{26} (1976), no.~1, 61--76},
}

@misc{ABMRW25,
    AUTHOR = {Affonso, Lucas and Bissacot, Rodrigo and Maia, Jo\~{a}o and
              Rodrigues, Jo\~{a}o F. and Welsch, Kelvyn},
     TITLE = {Cluster expansion and decay of correlations for
              multidimensional long-range {I}sing models},
      YEAR = {2025},
       URL = {https://arxiv.org/abs/2508.15666},
      NOTE = {Preprint, arXiv:2508.15666},
}

@book{Knuth_TAOCP1,
    AUTHOR = {Knuth, Donald E.},
     TITLE = {The art of computer programming. {V}ol. 1: {F}undamental algorithms},
   EDITION = {Third},
 PUBLISHER = {Addison-Wesley, Reading, MA},
      YEAR = {1997},
     PAGES = {xx+650},
      ISBN = {0-201-89683-4},
}

@article{GZ83,
    AUTHOR = {Greene, Curtis and Zaslavsky, Thomas},
     TITLE = {On the interpretation of {W}hitney numbers through arrangements of hyperplanes, zonotopes, non-{R}adon partitions, and orientations of graphs},
   JOURNAL = {Trans. Amer. Math. Soc.},
    VOLUME = {280},
      YEAR = {1983},
    NUMBER = {1},
     PAGES = {97--126},
  MRNUMBER = {712251},
       DOI = {10.1090/S0002-9947-1983-0712251-1},
}

@article{Tutte54,
    AUTHOR = {Tutte, W. T.},
     TITLE = {A contribution to the theory of chromatic polynomials},
   JOURNAL = {Canad. J. Math.},
    VOLUME = {6},
      YEAR = {1954},
     PAGES = {80--91},
       DOI = {10.4153/CJM-1954-010-9},
}

@book{Biggs_AGT,
    AUTHOR = {Biggs, Norman},
     TITLE = {Algebraic graph theory},
    SERIES = {Cambridge Mathematical Library},
   EDITION = {Second},
 PUBLISHER = {Cambridge University Press, Cambridge},
      YEAR = {1993},
     PAGES = {viii+205},
      ISBN = {0-521-45897-8},
   MRCLASS = {05C50},
  MRNUMBER = {1271140},
       DOI = {10.1017/CBO9780511608704},
}

@article{Griffiths69,
    AUTHOR = {Griffiths, Robert B.},
     TITLE = {Nonanalytic behavior above the critical point in a random {I}sing ferromagnet},
   JOURNAL = {Phys. Rev. Lett.},
    VOLUME = {23},
      YEAR = {1969},
     PAGES = {17--19},
       DOI = {10.1103/PhysRevLett.23.17},
}

\end{document}